\documentclass[12pt]{article}
\usepackage{amsmath}
\usepackage{amssymb}
\usepackage{amsthm}
\usepackage{amscd}
\usepackage{color}
\usepackage[mathscr]{eucal}

\DeclareMathOperator{\tr}{tr}
\DeclareMathOperator{\Hom}{Hom}

\DeclareMathOperator{\End}{End}
\DeclareMathOperator{\im}{Im}

\DeclareMathOperator{\id}{id}

\DeclareMathOperator{\slf}{S}

\newcommand{\C}{\mathbb{C}}
\newcommand{\Z}{\mathbb{Z}}

\newcommand{\CP}{\mathbb{P}}

\newcommand{\vac}{|\,0\,\rangle}

\newcommand{\frakh}{\mathfrak{h}}

\newcommand{\Res}[1]{\underset{#1=0}{\operatorname{Res}}\,}

\newcommand{\NO}{\,{\raise0.25em\hbox{$\mathop{\hphantom{\cdot}}%
\limits^{_{\circ}}_{^{\circ}}$}}\,}% The definition of the normal ordering.
\newcommand{\gch}{\operatorname{Ch}}

\theoremstyle{plain}
   \newtheorem{theorem}{Theorem}[subsection]%[subsection]
   \newtheorem{corollary}[theorem]{Corollary}
   \newtheorem{lemma}[theorem]{Lemma}
   \newtheorem{proposition}[theorem]{Proposition}
   \newtheorem{conjecture}[theorem]{Conjecture}
\theoremstyle{remark}
   \newtheorem{definition}[theorem]{Definition}
   \newtheorem{remark}[theorem]{Remark}

\numberwithin{equation}{section}

\date{\today}

\begin{document}
\begin{large}
\begin{center}
Some remarks on pseudo-trace functions for orbifold models \\
associated with  symplectic fermions
\end{center}
\end{large}
\vskip3ex
\begin{center}
Yusuke Arike and Kiyokazu Nagatomo\\
\vskip1ex
Department of Pure and Applied Mathematics\\
Graduate School of Information Science and  Technology\\
Osaka University\\
y-arike@cr.math.sci.osaka-u.ac.jp\\
\vskip1ex
Department of Pure and Applied Mathematics\\
Graduate School of Information Science and Technology\\
Osaka University\\
nagatomo@math.sci.osaka-u.ac.jp
\end{center}
\vskip2ex
\abstract{
We give a method to construct pseudo-trace functions for vertex operator algebras
satisfying Zhu's finiteness condition not through higher Zhu's algebras
and apply our method to the $\Z_2$-orbifold model associated with $d$-pairs
of symplectic fermions.
For $d=1$, we determine the dimension of the space of one-point functions.
For $d>1$, we construct $2^{2d-1}+3$ linearly independent one-point functions
and study their values at the vacuum vector.


\tableofcontents
\section{Introduction}\label{sect-1}
In this paper, we concern one-point functions associated with 
the vertex operator algebra $\mathscr{F}^+$ obtained as a
$\Z_2$-orbifold model of the symplectic fermionic vertex operator superalgebra with
$d$-pairs of symplectic fermions.
The vertex operator algebra $\mathscr{F}^+$ is intensively
studied by Abe \cite{A} and is known as the first example which is
not rational but satisfies Zhu's finiteness condition.

Let $V$ be a vertex operator algebra satisfying Zhu's finiteness condition.
In \cite{M}, it is shown that any one-point function is constructed in terms of 
a symmetric linear function on the $n$-th Zhu's algebra $A_n(V)$ 
introduced in \cite{DLM1} for sufficiently large integer $n$.
However, by the reason why it is not easy to determine the $n$-th Zhu's algebra $A_n(V)$, 
more simple way to obtain one-point functions is desired.

On the other hand, it is shown in \cite{MNT} that there exists a finite-dimensional
associative algebra whose category of finite-dimensional modules is equivalent
to the category of $V$-modules.
Then it is expected that the space of one-point functions is linearly isomorphic
to the one of symmetric linear functions on
the algebra. 
In fact, it is proved in \cite{T1} and \cite{T2} that
the space of one-point functions is isomorphic to the vector space of symmetric linear
functions on some algebra which is closely related  to conformal field theory.

We will take an another way, but adopting the spirit
proposed in \cite{M,T1,T2}. More precisely,
suppose we are given a $V$-module $M$. 
We consider a subalgebra $P$ of the endomorphism ring of $M$ such that $M$ is projective
as a left $P$-module.
This subalgebra $P$ is called a \textit{projective commutant}
in this paper. For any symmetric linear function $\varphi$ on $P$,
we are able to define a symmetric linear function $\varphi_M$
on the endomorphism ring $\End_P(M)$.
Then we can obtain a function $\varphi_M(J_0(a)q^{L_0-c_V/24})$ for any $a\in V$
which gives rise to a one-point function,
where $c_V$ denotes the central charge of $V$ and $J_0(a)$ denotes the zero-mode of $a\in V$. 
We remark that $L_0$ is not necessarily semisimple on $M$
so that this linear functional $\varphi_M$ may not be the trace on $M$.

The main interest in this paper is the space of one-point 
functions associated with the vertex operator algebra $\mathscr{F}^+$.
The vertex operator algebra $\mathscr{F}^+$ is an even part of the vertex operator superalgebra
$\mathscr{F}$ constructed from the $2d$-dimensional symplectic vector space $\frakh$.
The vertex operator algebra $\mathscr{F}^+$ for $d=1$ is known to be isomorphic 
to the triplet $\mathscr{W}(2)$-algebra.
It is shown in \cite{A} that there are four simple $\mathscr{F}^+$-modules
$\mathscr{F}^\pm$ and $\mathscr{F}_t^\pm$ for any integer $d\geq1$.
Moreover, two indecomposable $\mathscr{F}^+$-modules $\mathscr{F}^\pm_+$ are constructed
as the even and odd parts of the indecomposable $\mathscr{F}$-module $\mathscr{F}_+$.
We will find that the endomorphism ring $P=\End_{\mathscr{F}^+}(\mathscr{F}_+)$
is itself a projective commutant of $\mathscr{F}_+$
and show that $P$ is a symmetric algebra with a symmetric linear function $\varphi$.
It is well known that the center of a symmetric algebra is isomorphic to the
vector space of symmetric linear functions on the algebra.
The center of $P$ is explicitly determined and it is $(2^{2d-1}+1)$-dimensional,
which enables us to obtain a basis of the space of symmetric linear functions
on $P$.
Together with including ordinary trace functions defined on $\mathscr{F}^\pm_t$,
we have linearly independent $(2^{2d-1}+1+2)$ one-point functions.

For $d=1$, we show that one-point functions constructed here form
a basis of the space of one-point functions and that their values at the vacuum vector are also linearly independent.
In this proof, we use the fact that the dimension of the space of one-point functions 
is less than or equal to the dimension of the vector space of symmetric linear functions on Zhu's algebra $A_0(V)$.
Since Zhu's algebra $A_0(\mathscr{F}^+)$ with $d=1$ is known
(cf. \cite{A,AM2,NT2}), by using the explicit structure of $A_0(\mathscr{F}^+)$,
it is proved that the vector space of one-point functions is $5$-dimensional.

The notion of pseudo-characters is defined 
as values of one-point functions at the vacuum vector.
We will show that the space of one-point functions associated with 
$\mathscr{F}^+$ for $d=1$
is isomorphic to the space of pseudo-characters, while
in \cite{AM1}, it is proved that the space of pseudo-characters is $5$-dimensional
by showing that pseudo-characters are solutions of an ordinary differential equation.
For $d>1$, we will show that values at the vacuum vector of a part of one-point functions constructed 
in this paper are trivial.

This paper is organized as follows.
In Sect.~\ref{sect-2}, we recall definitions of modules for vertex operator algebras,
Zhu's finiteness condition and Zhu's algebras.
We also prove several properties for modules for vertex operator algebras
satisfying Zhu's finiteness condition.

In Sect.~\ref{sect-3}, we recall the definition of one-point functions
and prove that the dimension of the space of one-point functions is less than or equal to
the one of symmetric linear functions on Zhu's algebra $A_0(V)$.

In Sect.~\ref{sect-4}, we introduce a notion of pseudo-trace functions and prove that
pseudo-trace functions are one-point functions under Zhu's finiteness condition.

Sect.~\ref{sect-5} is devoted to the construction of the vertex operator algebra $\mathscr{F}^+$,
the classification of simple $\mathscr{F}^+$-modules,
and the construction of reducible indecomposable $\mathscr{F}^+$-modules.
The structure of Zhu's algebra described in \cite{A} are explained. 
It is worthy to mention that, for $d=1$, we have
\begin{equation}
A_0(\mathscr{W}(2))\cong A_0(\mathscr{F}^+)\cong \C\oplus M_2(\C)\oplus  M_2(\C)\oplus
\left\{
\begin{pmatrix}
a & b\\
0 & a
\end{pmatrix}
\,|\, a,\,b\in\C
\right\}.
\end{equation}

In Sect.~\ref{sect-6}, we apply our method given in Sect.~\ref{sect-4}
to construct pseudo-trace functions to $\mathscr{F}^+$-module 
$\mathscr{F}_+=\mathscr{F}_+^+\oplus\mathscr{F}_+^-$.
We explicitly determine the structure of the 
endomorphism ring $P$ of $\mathscr{F}_+$ and show that $\mathscr{F}_+$
is projective as a $P$-module. We take an appropriate $P$-coordinate system
of $\mathscr{F}_+$ (see \cite{Ar}) and construct pseudo-trace functions.

\section{Preliminaries}\label{sect-2}
In this section, we will recall definitions of modules for vertex operator algebras,
and Zhu's algebras.
Throughout this paper, $\Z$, $\Z_{\ge0}$, $\Z_{>0}$ and $\C$ denote the set of integers, non-negative integers,
positive integers and complex
numbers, respectively.

\subsection{Vertex operator algebras and their modules}\label{sect-2-1}
A vertex operator algebra is a quadruple $(V, Y, \vac, \omega)$
which consists of a $\Z$-graded vector space $V=\bigoplus_{n\in\Z}V_n$
over the field $\C$ of complex numbers,
a linear map
\begin{equation}
Y:V\to\End_\C(V)[[z,z^{-1}]]\quad(a\mapsto Y(a,z)=\sum_{n\in\Z}a_{(n)}z^{-n-1}),
\end{equation}
and two non-zero vectors $\vac\in V_0$ called the  {\it vacuum vector}  and $\omega\in V_2$
called the {\it Virasoro element} satisfying a number of conditions (see e.g. \cite{MaNa}).
For any $a\in V_n$, we write $|a|=n$.
We often refer $V$ to be a vertex operator algebra for short.
Throughout the paper,
we always assume that a vertex operator algebra $V$ is $\Z_{\ge0}$-graded, that is, $V=\bigoplus_{n=0}^\infty V_n$.

A {\it weak $V$-module} is a complex vector space $M$ equipped with a linear map
$Y_M:V\to\End_\C(M)[[z,z^{-1}]]$ which maps $a\in V$ to $Y_M(a,z)=\sum_{n\in\Z}a_{(n)}^Mz^{-n-1}$
satisfying the following:
\begin{enumerate}
\item[(1)]
For any pair $(a,u)\in V\times M$, there exists an integer $n_0$ such that
$a_{(n)}^Mu=0$ for all $n\ge n_0$.
\item[(2)]
$Y_M(\vac,z)=\id_M$.
\item[(3)]
For all $a, b\in V$ and $p,q,r\in\Z$,
the Borcherds identity holds:
\begin{equation}\label{eq-borcherds-voa}
\begin{split}
&\sum_{i=0}^\infty\binom{p}{i}(a_{(r+i)}b)_{(p+q-i)}^M\\
&\phantom{\sum_{i=0}^\infty(a_{(r+i)}}
=\sum_{i=0}^\infty(-1)^i\binom{r}{i}(a_{(p+r-i)}^Mb_{(q+i)}^M-(-1)^rb_{(q+r-i)}^Ma_{(p+i)}^M).
\end{split}
\end{equation}
\end{enumerate}

We will often omit the subscript and the superscript $M$ of $Y_M$ and $a_{(n)}^M\,(a\in V, n\in\Z)$
for simplicity.

Set $Y(\omega,z)=\sum_{n\in\Z}L_nz^{-n-2}$.
Then $\{L_n, \id_M\,|\,n\in\Z\}$ gives rise to a representation of the Virasoro algebra on $M$,
that is,
\begin{equation}\label{eq-Viro}
[L_m, L_n]=(m-n)L_{m+n}+\frac{m^3-m}{12}\delta_{m+n,0}c_V\id_M
\end{equation}
for all $m,n\in\Z$,
where $c_V\in\C$ is the central charge of $V$.
We also have
\begin{equation}\label{eq-derivative}
\frac{d}{dz}Y(a,z)=Y(L_{-1}a,z)
\end{equation}
for all $a\in V$ (see for instance \cite{MaNa}).

Setting $r=0$ and $p=0$ in the Borcherds identity, respectively, we have
\begin{align}
&[a_{(p)}, b_{(q)}]=\sum_{i=0}^\infty\binom{p}{i}(a_{(i)}b)_{(p+q-i)},\label{eq-commutator}\\
&(a_{(r)}b)_{(q)}=\sum_{i=0}^\infty(-1)^i\binom{r}{i}(a_{(r-i)}b_{(q+i)}
-(-1)^rb_{(q+r-i)}a_{(i)})\label{eq-associative}
\end{align}
called the \textit{commutator formula} and the \textit{associativity formula}.

For any homogeneous $a\in V$ and $n\in\Z$, we set $J_0(a)=a_{(|a|-1+n)}$ by $J_n(a)$
and extend it by linearity.
Note that $J_0(\vac)=\id_M$ and $J_n(\omega)=L_{n}$.
By the commutator formula and \eqref{eq-derivative}, we have
\begin{equation}\label{eq-gradation}
[L_0, J_n(a)]=-nJ_n(a)
\end{equation}
for all $a\in V$ and $n\in\Z$.

\begin{definition}[\cite{MNT}, cf. \cite{NT1}]\label{def-mods}
A module for a vertex operator algebra $V$ is a weak $V$-module $M$ satisfying
the following conditions:
\begin{enumerate}
\item[(1)]
$M$ is finitely generated as a weak $V$-module.
\item[(2)]
For any $u\in M$, there is an integer $m$ such that 
\begin{equation}
J_{n_1}(a^1)J_{n_2}(a^2)\dotsb J_{n_r}(a^r)u=0
\end{equation}
for all $a^1, a^2, \dotsc, a^r\in V$ and $n_1+n_2+\dotsb+n_r\ge m$.
\end{enumerate}
\end{definition}

We note that the integer $m$ appeared in Definition \ref{def-mods} (2)
must be positive since $J_0(\vac)=\id_M$.

\begin{remark}
In \cite{MNT}, a weak $V$-module satisfying Definition \ref{def-mods} (2)
is called an exhaustive module.
\end{remark}
\subsection{Zhu's finiteness condition}\label{sect-2-2}
In this subsection, we recall the definition of Zhu's finiteness condition introduced by Y.~Zhu (\cite{Z1})
and we give several properties of modules for a vertex operator algebra
satisfying Zhu's finiteness condition.

Let $V$ be a vertex operator algebra
and let $C_2(V)$ be the vector subspace of $V$ linearly spanned by elements
$a_{(-2)}b$ for all $a, b\in V$.
If $V/C_2(V)$ is finite-dimensional, we say that the vertex operator algebra $V$
satisfies {\it Zhu's finiteness condition}.

The following lemma is a modified version of \cite[Lemma 2.4]{M}
by using the operator $J_n(a)$.
The proof is similar to the one given in \cite{M}.

\begin{lemma}[cf. {\cite[Lemma 2.4]{M}, \cite[Lemma 8.3.1]{MNT}}]\label{lem-fermionic}
Let $V$ be a vertex operator algebra satisfying Zhu's finiteness condition
and let $M$ be a finitely generated weak $V$-module.
Then there is a finite-dimensional graded vector subspace $U$ of $V$
such that $V=C_2(V)+U$ and $M$ is linearly spanned by 
\begin{equation}\label{eq-fermionicp}
J_{n_1}(a^1)J_{n_2}(a^2)\dotsb J_{n_r}(a^r)u\,(r\in\Z_{\ge0}, a^i\in U, n_1<n_2<\dotsb<n_r, u\in B)
\end{equation}
where $B$ is a set of generators of $M$.
\end{lemma}

Under Zhu's finiteness condition, we can prove that any $V$-module
decomposes into a direct sum of finite-dimensional generalized eigenspaces for $L_0$.

\begin{proposition}[{\cite[Corollary 3.2.8]{NT1}}]\label{prop-mods-eigendecomp}
Let $V$ be a vertex operator algebra satisfying Zhu's finiteness condition
and let $M$ be a non-zero $V$-module.
Then there exist complex numbers $r_1, r_2, \dotsc, r_k$ with $r_i-r_j\not\in\Z$
for any $1\le i,j\le k$
and non-negative integers $d_1, d_2, \dotsc, d_k$
such that $M=\bigoplus_{i=1}^k\bigoplus_{n=0}^\infty M_{(r_i+n)}$
where $M_{(r_i+n)}=\{u\in M\,|\,(L_0-r_i-n)^{d_i+1}u=0\}$
and $M_{(r_i)}\neq0$ for all $1\le i\le k$.
Moreover each generalized eigenspace $M_{(r_i+n)}$ is finite-dimensional.
\end{proposition}
\begin{proof}
First we shall prove that $M$ decomposes into a direct sum of finite-dimensional generalized eigenspaces for $L_0$.

Let $B$ be a finite set of generators of $M$
and let $U$ be a finite-dimensional graded vector subspace of $M$
such that $V=C_2(V)+U$.
We may assume that $U$ contains the Virasoro element $\omega$.
By Lemma \ref{lem-fermionic} and the definition of $V$-modules, 
there exists a positive integer $m$ such that $M$ is linearly spanned by elements of the form
\eqref{eq-fermionicp} with $n_1+n_2+\dotsb+n_r<m$.

Let $W$ be the vector subspace of $M$ linearly spanned by
elements of the form \eqref{eq-fermionicp} with $0\le n_1+n_2+\dotsb+n_r<m$.
Since $U$ is finite-dimensional and $B$ is a finite set,
the vector space $W$ is finite-dimensional.
Since $W$ is $L_0$-invariant by \eqref{eq-gradation}, the vector space $W$ decomposes
into a direct sum of generalized eigenspaces for $L_0$.
We may assume that any $u\in B$ is a generalized eigenvector for $L_0$.
Then, for any $u\in B$, there exist a complex number $r_u$
and a non-negative integer $d_u$ such that $(L_0-r_u)^{d_u+1}u=0$.
By \eqref{eq-gradation}, we have
\begin{equation}
\begin{split}
&(L_0-r_u+n_1+n_2+\dotsb+n_r)^{d_u+1}J_{n_1}(a^1)J_{n_2}(a^2)\dotsb J_{n_r}(a^r)u\\
&=J_{n_1}(a^1)J_{n_2}(a^2)\dotsb J_{n_r}(a^r)(L_0-r_u)^{d_u+1}u\\
&=0
\end{split}
\end{equation}
for all $a^i\in U$, $u\in B$ and $n_1<n_2<\dotsb<n_r$,
which proves that $M$ decomposes into a direct sum of generalized eigenspaces for $L_0$.

Let $M_{(s)}$ be a generalized eigenspace for $L_0$ with eigenvalue $s$.
Then $M_{(s)}$ is linearly spanned by
$J_{n_1}(a^1)J_{n_2}(a^2)\dotsb J_{n_r}(a^r)u$
with $a^i\in U$, $u\in B$, $n_1<n_2<\dotsb<n_r$ and $n_1+n_2+\dotsb+n_r=r_u-s$.
This shows that each generalized eigenspace is finite-dimensional.

Set $P(B)=\{r_u\,|\,u\in B\}$.
Introduce the partial order on $\C$ as follows:
\begin{equation}\label{eq-partialc}
r\le s\text{ if and only if } s-r\in\Z_{\ge0}. 
\end{equation}
Then there exist minimal elements $r_1, r_2,\dotsc, r_k$ in $P(B)$
with respect to the partial order.
Thanks to Definition \ref{def-mods} (2), for any $r_i$, there exists a positive integer $m_i$
such that $M_{(r_i-m_i+n)}=0$ for all $n<0$ and $M_{(r_i-m_i)}\neq0$.
For any $1\le i\le k$, let $\{u_1^i, u_2^i, \dotsc, u_{\ell_i}^i\}$ be the subset of $B$ consisting of all elements
whose eigenvalues for $L_0$ are $r_i$.
Take the maximum element $d_i$ in the set $\{d_{u_1^i}, d_{u_2^i},\dotsc, d_{u_{\ell_i}^i}\}$.
Then we have $(L_0-r_i+m_i-n)^{d_i+1}M_{(r_i-m_i+n)}=0$
for all $n$.
This completes the proof.
%\qed
\end{proof}

Let $M=\bigoplus_{i=1}^k\bigoplus_{n=0}^\infty M_{(r_i+n)}\,(r_i-r_j\not\in\Z)$
be a module for a vertex operator algebra $V$ satisfying Zhu's finiteness condition.
Then, by \eqref{eq-gradation}, we have
$J_m(a)M_{(r_i+n)}\subset M_{(r_i+n-m)}$ for all $m\in\Z$ and $a\in V$.
Since $r_i-r_j\not\in\Z$ for any $1\le i,j\le k$,
it follows that $\bigoplus_{n=0}^\infty M_{(r_i+n)}$
is a $V$-submodule of $M$ for all $1\le i\le k$.

\begin{proposition}\label{prop-2-1-5}
Let $V$ be a vertex  operator algebra satisfying Zhu's finiteness condition
and let $M$ be a simple $V$-module.
Then there exists a complex number $r$ such that
$M=\bigoplus_{n=0}^\infty M_{r+n}$ with $M_{r}\neq0$ where 
$M_{r+n}$ is a finite-dimensional eigenspace for $L_0$ with eigenvalue $r+n$.
\end{proposition}

\subsection{Zhu's algebra $A_0(V)$}\label{sect-2-3}
In this subsection, we recall from \cite{Z1,DLM1} the definition of Zhu's algebra. 

Let $O(V)$ be the vector space linearly spanned by elements
\begin{equation}\label{3.1}
a\circ b
=\Res z Y(a,z)b\frac{(1+z)^{|a|}}{z^{2}}\,dz
\end{equation}
for any homogeneous $a\in V$ and $b\in V$, and set $A_0(V)=V/O(V)$.
We denote by $[a]$ the image of $a\in V$ in $A_0(V)$.
We define a bilinear operation $*$ on $V\times  V\rightarrow V$ 
by
\begin{equation}\label{2.10}
a*b=\Res z Y(a,z)b\frac{(1+z)^{|a|}}{z}\,dz
\end{equation}
for homogeneous $a\in V$ and $b\in V$.

\begin{theorem}[\cite{Z1}]\label{prop-1-2}
The bilinear operation $*$ induces a structure of an associative $\C$-algebra 
on $A_0(V)$.
Moreover, $[\vac]$ is the unity of $A_0(V)$
and $[\omega]$ is a central element.
\end{theorem}

\begin{theorem}[\cite{Z1}]
\label{theorem2.2.3}
Let $V$ be a vertex operator algebra and let $M$ be a $V$-module.
\begin{enumerate}
\item[{\rm (1)}]
The linear map $J_0:V\rightarrow\End_\C(\Omega(M))$ defined by $a\mapsto J_0(a)$ induces a representation of 
$A_0(V)$ on $\Omega(M)$ where 
$\Omega(M)
=\{m\in M\,|\, J_n(a)m=0\text{ for any } a\in V\text{ and integers } n>0\}$.
\item[{\rm (2)}]
The operation $\Omega$ induces a bijection between the complete set of inequivalent 
simple $V$-modules and the set of inequivalent simple $A_0(V)$-modules.
\end{enumerate}
\end{theorem}

\section{The space of one-point functions}\label{sect-3}
In this section, we recall the definition of one-point functions
and discuss their properties.

\subsection{One-point functions on elliptic curves}\label{sect-3-1}
In this subsection, we shall recall the definition 
of one-point functions on elliptic curves.

The Eisenstein series $G_{2k}(\tau)\,(k\ge1)$ are series defined by
\begin{equation}
G_{2k}(\tau)=\sum_{(m,n)\neq(0,0)}\frac{1}{(m\tau+n)^{2k}}
\end{equation}
for $k\ge2$
and
\begin{equation}
G_2(\tau)=\frac{\pi^2}{3}+\sum_{m\in\Z-\{0\}}\sum_{n\in\Z}\frac{1}{(m\tau+n)^2}.
\end{equation}
It is well-known that the Eisenstein series have the $q$-expansions
\begin{equation}
G_{2k}(\tau)=2\zeta(2k)+\frac{2(2\pi i\tau)^{2k}}{(2k-1)!}\sum_{n=1}^\infty
\sigma_{2k-1}(n)q^n
\end{equation}
for all $k\ge1$,
where $\zeta(s)$ is the Riemann zeta function, $\sigma_{k}(n)=\sum_{d|n}d^k$,
$q=e^{2\pi i\tau}$ and $\tau\in\mathscr{H}=\{\tau\in\C\,|\,\im\tau>0\}$.

For any elliptic curve $E_\tau=\C/(\Z\oplus\Z\tau)$ with a modulus $\tau$ in $\mathscr{H}$,
we have two canonical coordinate systems:
the one is the coordinate $z$ of $\C$ and the other is the coordinate 
$e^{2\pi iz}$ of $\CP^1$. There exist two vertex operator algebra 
structures corresponding to these two coordinate systems, respectively
(cf. \cite{FB} for instance).

\begin{theorem}[{\cite[Theorem 4.2.1]{Z1}}]\label{thm-1}
Let $(V,Y,\vac,\omega)$ be a vertex operator algebra of central charge $c_V$. 
Define $Y[a, z] = Y(a,\,e^{2\pi iz}-1)e^{2\pi i|a|z}$
for every homogeneous $a\in V$ and set 
$\tilde{\omega} =(2\pi i)^2\left(\omega - \frac{c_V}{24}\vac\right)$.
Then the quadruple $(V, Y[\,], \vac, \tilde{\omega})$ is 
a vertex operator algebra.
\end{theorem}

For any $a\in V$, we denote $Y[a,z]=\sum_{n\in\Z}a_{[n]}z^{-n-1}$ and
set $L_{[n]}=\tilde{\omega}_{[n+1]}$ for any integer $n$.
Then the set
\begin{equation}
V_{[n]}=\{v\in V\,|\,L_{[0]}a=na\}.
\end{equation}
is not equal to $V_n$ in general,
however, it is shown in \cite{DLM3}
that for each non-negative integer $n$ we have
\begin{equation}\label{eq-negativestates}
\bigoplus_{k\leq n}V_k=\bigoplus_{k\leq n}V_{[k]}.
\end{equation}
Since $V$ is assumed to be $\Z_{\ge0}$-graded,
it follows from \eqref{eq-negativestates} that
the vertex operator algebra $(V, Y[\,], \vac, \tilde{\omega})$
is also $\Z_{\ge0}$-graded.

Let $O_q(V)$ be the $\mathbb{C}[G_4, G_6]$-submodule of
$V\otimes \C[G_4, G_6]$ generated by
\begin{align}
&a_{[0]}b,\label{e-2-2}\\
&a_{[-2]}b+\sum_{k = 2}^{\infty}(2k-1)
a_{[2k-2]}b\otimes G_{2k}(\tau)\label{e-2-1}
\end{align}
for all homogeneous $a\in V$ and $b\in V$.
Then we can define one-pint functions.

\begin{definition}
Let $\Gamma\subset\C$ be an integral lattice
and let $E_\Gamma=\C/\Gamma$ be the corresponding elliptic curve.
A map $S:V\otimes \C[G_4,\,G_6]\times\mathscr{H}\rightarrow\C$
satisfying the following conditions is called 
a \textit{one-point function on the elliptic curve $E_\Gamma$}.
\begin{enumerate}
\item[(1)]
For any $a\in V\otimes \C[G_4,G_6]$, the function $S(a,\tau)$ is
holomorphic in $\tau\in\mathscr{H}$.
\item[(2)]
$S(\sum_{i}a_i\otimes f_i(\tau),\tau)
=\sum_{i}f_i(\tau)S(a_i,\tau)$ for all $a_i\in V$ and 
$f_i\in\mathbb{C}[G_4,\,G_6]$.
\item[(3)]
$S(a,\tau) = 0$ for all $a\in O_q(V)$.
\item[(4)]
For any $a\in V_{[n]}$, the following relation holds:
\begin{equation}\label{e-2-3}
S (L_{[-2]}a, \tau)
=(2\pi i )^2q\frac{d}{dq}S(a,\tau)
+\sum_{k = 1}^{\infty}G_{2k}(\tau)S(L_{[2k-2]}a,\tau)
\end{equation}
where $q=e^{2\pi i\tau}$.
\end{enumerate}
\end{definition}
We denote the space of one-point functions by $\mathscr{C}(V)$.

Let $V$ be a vertex operator algebra satisfying Zhu's finiteness condition
and let $S$ be a non-zero one-point function.
It is shown in \cite[pp.~85]{M} (also see \cite[pp.~295]{Z1} and \cite[Theorem 6.5]{DLM3})
that there exist complex numbers $r_1, r_2, \dotsc, r_d$ with $r_i-r_j\not\in\Z$ for $i\neq j$
and non-negative integers $N_1, N_2, \dotsc, N_d$
such that
\begin{align}
&S(a,\tau)=\sum_{i=1}^{d}S_i(a,\tau)q^{r_i},\\
&S_i(a,\tau)=\sum_{j=0}^{N_i}S_{ij}(a,\tau)(2\pi i\tau)^j, \label{eq-reduced-1}\\
&S_{ij}(a,\tau)=\sum_{k=0}^\infty S_{ijk}(a)q^k
\end{align}
for all $a \in V$.
It is obvious that $S_{ijk}\in\Hom_\C(V,\C)$.

\begin{remark}\label{remark-1019}
Since $r_i-r_j\not\in\Z$ for $i\neq j$,
the map $V\otimes\C[G_4, G_6]\times\mathscr{H}\to\C$
defined by $(a\otimes f, \tau)\mapsto S_i(a\otimes f,\tau)q^{r_i}$
is a one-point function (cf. \cite{M}).
\end{remark}

Let $S$ be a one-point function.
If there is a complex number $r$ such that
\begin{equation}
S(a,\tau)=\sum_{j=0}^d\sum_{k=0}^\infty S_{jk}(a)q^{r-c_V/24+k}(2\pi i\tau)^j
\end{equation}
for all $a\in V$ with $S_{00}\neq0$ where $c_V$ denotes the central charge of $V$,
then we call $S$ a {\it one-point function of conformal weight $r$}.

\begin{definition}
We set $\gch(V)=\{S(\vac,\tau)\,|\,S\in\mathscr{C}(V)\}$.
We call any element of $\gch(V)$ a
\textit{pseudo-character} (It is called a generalized character in \cite{M}).
\end{definition}

Note that, by definition, there is a surjective linear map
$\mathscr{C}(V)\to\gch(V)$
defined by $S\mapsto S(\vac,\tau)$.
If $V$ satisfies Zhu's finiteness condition, then the vector space $\gch(V)$ contains
characters of simple $V$-modules and it is finite-dimensional 
(see \cite{M,Z1}).

\subsection{Symmetric linear functions on associative algebras}\label{sect-3-2}
In order to state the relationship between one-pint functions and Zhu's algebra,
we recall the notion of symmetric linear functions on associative algebras.

Let $A$ be a finite-dimensional associative $\C$-algebra
and let $Z(A)$ be the center of $A$.
A linear function $\phi:A\to\C$ is called a {\it symmetric linear function} if
$\phi(ab)=\phi(ba)$ for all $a, b\in A$.
We denote the vector space of symmetric linear functions on $A$ by $\slf^A$.
The vector space $\slf^A$ is canonically 
a $Z(A)$-module by the action
\begin{equation}\label{eq-echoes-act3}
(c\cdot\varphi)(a)=\varphi(ca)\text{ for all } a\in A.
\end{equation}

\begin{proposition}\label{prop:dimension-1}
Let $A$ be a finite-dimensional associative $\C$-algebra 
and let $\omega$ be a non-zero central element.
Then there exist complex numbers $r_1, r_2,\dotsc, r_k$
and non-negative integers $d_1, d_2, \dotsc, d_k$
such that
$\prod_{i=1}^k(\omega -r_i)^{d_i+1}=0$. 
In particular, $A=\bigoplus_{i=1}^kA_{r_i}$
where $A_{r_i}=\{a\in A\,|\,(\omega-r_i)^{d_i+1}a=0\}$.
\end{proposition}
\begin{proof}
Let $M$ be a simple left $A$-module. 
Then there exists a primitive idempotent 
$e$ such that $M\cong Ae/J(A)e$
where $J(A)$ is the Jacobson radical of $A$.
Since $\omega$ acts on $Ae/J(A)e$ 
as a scalar $r$ by Schur's lemma,
we see that $(\omega-r)e\in J(A)e\subset J(A)$.
Because $J(A)$ is nilpotent and $\omega\in Z(A)$, there exists a non-negative integer
$d$ such that $(\omega -r)^{d+1}e=0$.
Therefore $(\omega -r)^{d+1}Ae=0$.
Since $A$ is a direct sum of left $A$-modules generated by primitive idempotents,
we have the assertion.
%\qed
\end{proof}

By \eqref{eq-echoes-act3} and Proposition \ref{prop:dimension-1}, we have:

\begin{proposition}\label{prop:dimension-2}
Let $A$ be a finite-dimensional associative $\C$-algebra and
let $\omega$ be a non-zero central element of $A$.
Then there exist complex numbers $r_1, r_2,\dotsc, r_k$
and non-negative integers $d_1, d_2,\dotsc,d_k$
such that $\slf^A=\bigoplus_{i=1}^k\slf^A_{r_i}$ where 
$(\omega-r_i)^{d_i+1}\slf^A_{r_i}=0$.
\end{proposition}

\subsection{Zhu's algebras and one-point functions}\label{sect-3-3}
Let $V$ be a vertex operator algebra satisfying Zhu's finiteness condition.
It is shown in \cite[Proposition 3.6]{DLM3} that $A_0(V)$ is finite-dimensional.
Then, by Theorem \ref{theorem2.2.3}, there are finitely many simple $V$-modules.
Let $\{M^1, M^2, \dotsc, M^k\}$ be the complete set of all inequivalent simple $V$-modules.
By Proposition \ref{prop-2-1-5}, for any $1\le i\le k$, there is a complex number $r_i$
such that $M^i=\bigoplus_{n=0}^\infty M_{r_i+n}^i$ with $M^i_{r_i}\neq0$
where each $M_{r_i+n}^i$ is a finite-dimensional eigenspace for $L_0$ with eigenvalue $r_i+n$.
The complex number $r_i$ is called the {\it conformal weight} of $M^i$.
Let $\Lambda$ be the subset of $\C$ consisting of  conformal weights of all simple $V$-modules.
Note that the cardinality of $\Lambda$ is equal to or smaller than $k$.

\begin{remark}
We use the terminology ``conformal weight'' for simple $V$-modules
and one-point functions.
By the result given in \cite[Lemma 5.7]{M} (see Lemma \ref{lem:slf} in this paper),
we will find that, if there is a one-point function of conformal weight $r$,
then $r\in\Lambda$.
\end{remark}

Throughout the paper, we denote the space of symmetric linear functions on $A_0(V)$ by $\slf^V$.
Then, by Proposition \ref{prop:dimension-2}, we obtain
\begin{equation}
\slf^V=\bigoplus_{r\in\Lambda}\slf^V_{r}
\end{equation}
where 
\begin{equation}
\slf^V_{r}=
\{\phi\in\slf^V\,|\,\phi(([\omega]-r)^{d_r+1}\ast 
[a])=0,\,\forall [a]\in A_0(V)\}.
\end{equation}

\begin{lemma}[{\cite[pp.82]{M}}]\label{lem:red}
Let $S$ be a one-point function.
Suppose that $S(a,\tau)$ is expressed as
$S(a,\tau)=\sum_{j=0}^d\sum_{k=0}^{\infty}S_{jk}(a)q^{r+k}(2\pi i\tau)^j$
for all $a\in V$ where $r\in\C$.
If $S_{00}=0$, then $S_{j0}=0$ for all $0\leq j\leq d$.
\end{lemma}

\begin{lemma}[\cite{Z1}, {\cite[Lemma 5.7]{M}}]\label{lem:slf}
Let $S$ be a one-point function.
Suppose that $S(a,\tau)$ is expressed as
$S(a,\tau)=\sum_{j=0}^d\sum_{k=0}^{\infty}S_{jk}(a)q^{r-c_V/24+k}
(2\pi i\tau)^j$ for all $a\in V$ where $r\in\C$.
Then $S_{00}\in\slf^V$ and $S_{00}((\omega-r)^{d+1}\ast a)=0$ for all $a\in V$.
In particular, if $S_{00}\neq0$, then $r\in\Lambda$ and $S_{00}$ belongs to $\slf^V_{r}$.
\end{lemma}

The following theorem is one of the main results in \cite{M}.

\begin{theorem}[{\cite[Theorem 5.5]{M}}]\label{thm:pseudo-trace}
Let $V$ be a vertex operator algebra satisfying Zhu's finiteness condition
and let $c_V$ be the central charge of $V$.
Suppose that any simple $V$-module is infinite-dimensional.
Then the vector space $\mathscr{C}(V)$ has a basis 
$\{S^{r,i_r}\,|\,r\in\Lambda,\, 1\leq i_r\leq k_r\}$ where
\begin{equation}
S^{r,i_r}(a,\tau)
=\sum_{j=0}^{d_{i_r}}\sum_{k=0}^\infty S_{jk}^{r,i_r}(a)q^{r-c_V/24+k}(2\pi i\tau)^j
\end{equation}
for all $a\in V$
with $S^{r,\,i_r}_{00}\not=0$.
Moreover, any one-point function of conformal weight $r_1$
is a linear combination of $S^{r_2, i_{r_2}}$ such that 
$\operatorname{Re}(r_2)\geq\operatorname{Re}(r_1)$.
\end{theorem}

\begin{remark}
The condition ``any simple $V$-module is infinite-dimensional'' in Theorem \ref{thm:pseudo-trace}
is not explicitly assumed in \cite{M}.
The basis of $\mathscr{C}(V)$ in Theorem \ref{thm:pseudo-trace}
consists of pseudo-trace functions in the sense of \cite{M}.
In order to define pseudo-trace functions in \cite{M},
the integer $N$ such that $M_{r_i+n}^i\neq0$ for all $n>N$ and $1\le i\le k$
plays an important role, where $\{M^1, M^2,\dotsc, M^k\}$ 
is the complete set of all inequivalent simple $V$-modules (see \cite[pp.~71]{M}).
However, if there is a finite-dimensional simple $V$-module,
we cannot choose such an integer 
and therefore may not obtain the basis as in Theorem \ref{thm:pseudo-trace}.
For example, the vertex operator algebra $\mathscr{W}_{2,3}$
satisfies Zhu's finiteness condition but admits a finite-dimensional simple module (see \cite{AM3}).
\end{remark}

By Theorem \ref{thm:pseudo-trace}, we obtain an upper bound of the dimension of the space of one-point functions:

\begin{theorem}\label{theorem3.3.4}
Let $V$ be a vertex operator algebra satisfying Zhu's finiteness condition.
Suppose that any simple $V$-module is infinite-dimensional.
Then $\dim_\C\mathscr{C}(V)\leq\dim_\C\slf^V$.
\end{theorem}
\begin{proof}
By Theorem \ref{thm:pseudo-trace}, we can choose a basis 
$\{S^{r,\,i_r}\,|\,r\in\Lambda,\,1\leq i_r\leq k_r\}$ of $\mathscr{C}(V)$ such that
\begin{align}
&S^{r,\,i_r}(a,\tau)
=\sum_{j=0}^{d_{i_r}}\sum_{k=0}^\infty S^{r,\,i_r}_{jk}(a)q^{r-c_V/24+k}(2\pi i\tau)^j,\\
&S^{r,\,i_r}_{jk}\in\Hom_\C(V,\C),\qquad S^{r,\,i_r}_{00}\not=0.
\end{align}

Suppose that $\dim_\C\mathscr{C}(V)>\dim_\C\slf^V$.
Then there exists a conformal weight $r\in\Lambda$ such that
$\dim_\C\slf^V_{r}<k_r$ since
$\dim_\C\mathscr{C}(V)=\sum_{r\in\Lambda}k_r$.
By Lemma \ref{lem:slf}, we see that the set
$\{S_{00}^{r,\,i_r}\,|\,1\leq i_r\leq k_r\}$ is contained in $\slf^V_{r}$.
Therefore the set $\{S_{00}^{r,\,i_r}\,|\,1\leq i_r\leq k_r\}$ 
is not linearly independent so that there exist  $i$ and complex numbers 
$a_{j_r}$ such that $S_{00}^{r,\,i}=\sum_{j_r\neq i}a_{j_r}S_{00}^{r,\,j_r}$. Set 
\begin{equation}\label{3.4}
T=S^{r,\,i}-\sum_{j_r\neq i}a_{j_r}S^{r,\,j_r}.
\end{equation}
Hence we can express $T$ as
\begin{equation}
T(a,\tau)=\sum_{j=0}^d\sum_{k=0}^\infty T_{jk}(a)q^{r-c_V/24+k}(2\pi i\tau)^j
\end{equation}
for all $a\in V$.
Then $T_{00}=0$ so that 
we have $T_{j0}=0$ 
for all $0\leq j\leq d$ by Lemma \ref{lem:red}. 
Therefore the one-point function 
$T$ is rewritten as
\begin{equation}
T(a,\tau)
=\sum_{j=0}^d\sum_{k=0}^\infty T_{jk}(a)q^{s-c_V/24+k}(2\pi i\tau)^j,
\end{equation}
for all $a\in V$
where $\operatorname{Re}(s)>\operatorname{Re}(r)$. 
By Theorem \ref{thm:pseudo-trace}, $T$ is a linear combination of one-point functions in the basis
whose conformal weights have real parts greater than or equal to $\operatorname{Re}(s)$.
Then,
by \eqref{3.4}, we see that $S^{r,\,i}$ is a linear combination of
other elements of the basis of $\mathscr{C}(V)$.
Therefore the set $\{S^{r,\,i_r}\,|\,r\in\Lambda,\,1\leq i\leq k_r\}$ is not linearly independent. 
This is a contradiction.
%\qed
\end{proof}

\section{Pseudo-trace functions}\label{sect-4}
In this section, we introduce a generalization of a notion of pseudo-trace functions
defined in \cite{M}.
For this purpose, we recall pseudo-trace maps and introduce a notion of
a projective commutant of modules for a vertex operator algebra $V$.
Then we prove that
any pseudo-trace function is a one-point function if $V$ satisfies Zhu's finiteness condition.

\subsection{Pseudo-trace maps}\label{sect-4-1}
In this subsection, we will recall the notion of pseudo-trace maps (see e.g. \cite{Ar}). 

Let $A$ be a finite-dimensional associative $\C$-algebra and
Let $W$ be a finitely generated projective left $A$-module.
It is well known that there exists a pair of sets 
$\{u_i \}_{i=1}^{n}\subset W$ and 
$\{ f_i \}_{i=1}^{n}\subset\Hom_A(W, A)$ such that
$w = \sum_{i=1}^{n}f_i(w)u_i $ for all $w\in W$. 
The set $\{u_i, f_i\}_{i=1}^{n}$ is called an 
\textit{$A$-coordinate system} of $W$ (cf. \cite[Chapter II, \S\,2.6, Proposition 12]{Bour}).

Let $\{u_i,f_i\}_{i=1}^n$ be an $A$-coordinate system of $W$.
For any $\phi\in\slf^A$, we define the linear map
$\phi_W:\End_A(W)\to\C$ by $\phi_{W} (\alpha)=\sum_{i = 1}^{n}\phi(f_i(\alpha(u_i)))$.
The map $\phi_{W}$ is called the \textit{pseudo-trace map}.

\begin{proposition}[{cf.\,\cite[\S\,2]{Ar}}]
Let $A$ be a finite-dimensional associative $\C$-algebra,
$W$ a finitely generated projective left
$A$-module and $\phi$ a symmetric linear function on $A$.  
Then the pseudo-trace map $\phi_{W}$ is independent of 
the choice of $A$-coordinate systems.
\end{proposition}

The following proposition will be needed
to prove that pseudo-trace functions are one-point functions.

\begin{proposition}[cf.\,\cite{Brou}]\label{proposition:2.2.5}
Let $A$ be a finite-dimensional associative $\C$-algebra
and let $\phi$ be a symmetric linear function on 
$A$. Suppose $P$ and $Q$ are finitely generated projective
left $A$-modules. Then we have 
$\phi_{P}(\beta\circ \alpha)=\phi_{Q}(\alpha\circ\beta)$
for all $\alpha\in\Hom_{A}(P,Q)$ and $\beta\in\Hom_{A}(Q,P)$.
In particular, $\phi_P$ is a symmetric linear function on $\End_A(P)$.
\end{proposition}
\begin{proof}
Let $\{u_i, f_i\}_{i=1}^m$ and $\{v_j, g_j\}_{j=1}^n$
be $A$-coordinate systems of $P$ and $Q$, respectively.
Then we have
\begin{equation}\label{eq-expandmap}
\alpha(u_i)=\sum_{j=1}^ng_j(\alpha(u_i))v_j,\,
\beta(v_j)=\sum_{i=1}^mf_i(\beta(v_j))u_i,
\end{equation}
It follows from \eqref{eq-expandmap} that
\begin{equation}
\begin{split}
\phi_{P}(\beta\circ\alpha)&=\sum_{i=1}^m\phi(f_i\circ\beta\circ\alpha(u_i))
=\sum_{i=1}^m\sum_{j=1}^n\phi(f_i\circ\beta(g_j(\alpha(u_i))v_j))\\
&=\sum_{i=1}^m\sum_{j=1}^n\phi(g_j(\alpha(u_i))f_i(\beta(v_j))
=\sum_{i=1}^m\sum_{j=1}^n\phi(f_i(\beta(v_j)g_j(\alpha(u_i)))\\
&=\sum_{i=1}^m\sum_{j=1}^n\phi(g_j\circ\alpha(f_i(\beta(v_j))u_i))
=\sum_{j=1}^n\phi(g_j\circ\alpha\circ\beta(v_j))\\
&=\phi_Q(\alpha\circ\beta),
\end{split}
\end{equation}
which shows the proposition.
%\qed
\end{proof}
\subsection{Projective commutants}\label{sect-4-2}

Let $M$ be a module for a vertex operator algebra $V$.
A subalgebra $P$ of $\End_V(M)$ is called a {\it projective commutant of $M$}
if $M$ is projective as a left $P$-module.

Suppose that $V$ satisfies Zhu's finiteness condition.
By Proposition \ref{prop-mods-eigendecomp},
any $V$-module $M$ decomposes into a direct sum of generalized eigenspaces for $L_0$
and each generalized eigenspace is finite-dimensional.
It is shown in \cite[Proposition 5.9.1]{NT1} that $\End_V(M)$ is finite-dimensional,
which shows that any projective commutant $P$ of $M$ is finite-dimensional.
Since any $\alpha\in\End_V(M)$ preserves generalized eigenspaces for $L_0$,
we have the following:

\begin{proposition}\label{prop-echoes-act5}
Let $V$ be a vertex operator algebra satisfying Zhu's finiteness condition
and let $M$ be a $V$-module.
Then any projective commutant $P$ of $M$ is finite-dimensional
and each generalized eigenspace for $L_0$ is a finite-dimensional
projective left $P$-module.
\end{proposition}

\subsection{Pseudo-trace functions and one-point functions}\label{sect-4-3}
Let $V$ be a vertex operator algebra satisfying Zhu's finiteness condition
and let $M$ be a $V$-module.
Then we can assume that
there exist a complex number $r$ and non-negative integer $d$ such that
$M=\bigoplus_{n=0}^\infty M_{(r+n)}$ with $M_{(r)}\neq0$ where
$(L_0-r-n)^{d+1}M_{(r+n)}=0$ and $\dim_\C M_{(r+n)}<\infty$ for all $n\ge0$ (see Sect.~\ref{sect-2-2}).
We define the operator $q^{L_0}$ on $M$ by
\begin{equation}\label{eq-qoperator}
q^{L_0}=\sum_{j=0}^d\frac{1}{j!}(L_0-r-n)^jq^{r+n}(2\pi i\tau)^j \text{ on } M_{(r+n)}
\end{equation}
where $q=e^{2\pi i\tau}$ and $\tau\in\mathscr{H}$.

Let $P$ be a projective commutant of $M$
and let $\phi$ be a symmetric linear function on $P$.
Since each $M_{(r+n)}$ is a finite-dimensional projective left $P$-module by Proposition \ref{prop-echoes-act5},
we can define the pseudo-trace map $\phi_{M_{(r+n)}}:\End_P(M_{(r+n)})\to\C$
for all $n\ge0$.
Note that
$J_m(a)\in\Hom_P(M_{(r+n)}, M_{(r+n-m)})$
for all $a\in V$, $m\in\Z$ and non-negative integers $n$.
We can set
\begin{equation}
\phi_{M_{(r+n)}}(J_0(a)q^{L_0})=\sum_{j=0}^d\frac{1}{j!}\phi_{M_{(r+n)}}(J_0(a)(L_0-r-n)^j)q^{r+n}
(2\pi i\tau)^j
\end{equation}
for all $a\in V$ and $n\ge0$,
\begin{equation}\label{eq-20111020-1}
\phi_M Y(a,z)q^{L_0}=z^{-|a|}\sum_{n=0}^{\infty}\phi_{M_{(r+n)}}\left(J_0(a)q^{L_0}\right)
\end{equation}
and
\begin{equation}\label{eq-20111020-2}
\begin{split}
&\phi_M Y(a,z)Y(b,w)q^{L_0}\\
&\qquad\qquad
=\sum_{n=0}^{\infty}\sum_{m\in\Z}^\infty z^{m-|a|}w^{-m-|b|}\phi_{M_{(r+n)}}\left(J_{-m}(a)J_{m}(b)q^{L_0}\right)
\end{split}
\end{equation}
for homogeneous $a, b\in V$.
Let $\{u_i^n, \alpha_i^n\,|\,1\le i\le k_n\}$ be a $P$-coordinate system of $M_{(r+n)}$.
We can extend $\alpha_i^n$ to the $P$-homomorphism $M\to P$ by letting $\alpha_i(M_{(r+m)})=0$
if $n\neq m$.
Then \eqref{eq-20111020-1} and \eqref{eq-20111020-2} respectively become
\begin{align}
&\phi_M Y(a,z)q^{L_0}=\sum_{n=0}^{\infty}
\phi_{M_{(r+n)}}(Y(a,z)q^{L_0})\label{eq-0405-1}\\
\intertext{and}
&\phi_M Y(a,z)Y(b,w)q^{L_0}=\sum_{n=0}^{\infty}
\phi_{M_{(r+n)}}(Y(a,z)Y(b,w)q^{L_0})\label{eq-0405-2}
\end{align}
for any homogeneous $a, b\in V$.

\begin{proposition}\label{prop-cyclic}
Let $M=\bigoplus_{n=0}^\infty M_{(r+n)}$ be a $V$-module
and let $P$ be a projective commutant of $M$.
Then we have
\begin{equation}
\begin{split}
\phi_M Y(a,z)Y(b,w)q^{L_0}
&=q^{|a|}\phi_MY(b,w)Y(a,zq)q^{L_0}\\
&=q^{-|b|}\phi_M Y(b,wq^{-1})Y(a,z)q^{L_0}
\end{split}
\end{equation}
for any homogeneous $a, b\in V$, $m\in\Z$ and $\phi\in\slf^P$.
\end{proposition}
\begin{proof}
Since $M_{(r+n)}=0$ for $n<0$, we have
$\phi_{M_{(r+n)}}(J_{-m}(a)J_{m}(b)q^{L_0})=0\,(0\leq n\leq m-1)$
for any positive integer $m$.
By Proposition \ref{proposition:2.2.5} and the fact that
\begin{equation}
\begin{split}
J_m(b)q^{L_0}&=\sum_{j=0}^{d_r}\frac{1}{j!}J_m(b)(L_0-r-n)^jq^{r+n}(2\pi i\tau)^j\\
&=\sum_{j=0}^{d_r}\frac{1}{j!}(L_0-r+m-n)^jJ_m(b)q^{r+n}(2\pi i\tau)^j\\
&=q^mq^{L_0}J_m(b)
\end{split}
\end{equation}
on each $M_{(r+n)}$, we obtain
\begin{equation}\label{3.20}
\phi_{M_{(r+n)}}(J_{-m}(a)J_{m}(b)q^{L_0})
=q^{m}\phi_{M_{(r+n-m)}}(J_{m}(b)J_{-m}(a)q^{L_0}).
\end{equation}

Remark that the left-hand side of \eqref{3.20} is zero for any 
$0\leq n\leq m-1\,(m>0)$ and so is the right-hand side. 
Therefore we get
\begin{equation}
\begin{split}
\sum_{n=0}^\infty\phi_{M_{(r+n)}}(J_{-m}(a)J_{m}(b)q^{L_0})
&=\sum_{n=m}^\infty\phi_{M_{(r+n)}}(J_{-m}(a)J_{m}(b)q^{L_0})\\
&=\sum_{n=m}^\infty q^{m}\phi_{M_{(r+n-m)}}(J_{m}(b)J_{-m}(a)q^{L_0})\\
&=\sum_{n=0}^\infty q^{m}\phi_{M_{(r+n)}}(J_{m}(b)J_{-m}(a)q^{L_0})\\
\end{split}
\end{equation}
for any non-negative integer $m$.
Since
\[
J_{-m}(a)M_{(r+n)}\subseteq M_{(r+n+m)}=0\quad(0\leq n\leq-m-1)
\]
for any negative integer $m$, it follows that
\begin{equation}
\phi_{M_{(r+n)}}(J_{m}(b)J_{-m}(a)q^{L_0})=0.
\end{equation}
Hence we obtain 
\begin{equation}
\begin{split}
\sum_{n=0}^\infty\phi_{M_{(r+n)}}(J_{-m}(a)J_{m}(b)q^{L_0})
&=\sum_{n=0}^\infty q^{m}\phi_{M_{(r+n-m)}}(J_{m}(b)J_{-m}(a)q^{L_0})\\
&=\sum_{n=-m}^\infty q^{m}\phi_{M_{(r+n)}}(J_{m}(b)J_{-m}(a)q^{L_0})\\
&=\sum_{n=0}^\infty q^{m}\phi_{M_{(r+n)}}(J_{m}(b)J_{-m}(a)q^{L_0})
\end{split}
\end{equation}
for any negative integer $m$.
As a consequence, we have
\begin{equation}\label{6.4}
\begin{split}
&\phi_M Y(a,z)Y(b,w)q^{L_0}\\
&\qquad=\sum_{m\in\Z}
\sum_{n=0}^\infty\phi_{M_{(r+n)}}(J_{-m}(a)J_{m}(b)q^{L_0})
z^{m-|a|}w^{-m-|b|}\\
&\qquad=\sum_{m\in\Z}\sum_{n=0}^\infty q^{m}\phi_{M_{(r+n)}}(J_{m}(b)J_{-m}(a)q^{L_0})
z^{m-|a|}w^{-m-|b|}\\
&\qquad=q^{|a|}\phi_MY(b,w)Y(a,zq)q^{L_0}
\end{split}
\end{equation}
and 
\begin{equation}
\begin{split}
&\phi_M Y(a,z)Y(b,w)q^{L_0}\\
&\phantom{\phi_M Y(u,z)}=\sum_{m\in\Z}\sum_{\ell=0}^\infty q^{m}\phi_{M_{(r+n)}}(J_{m}(b)J_{-m}(a)
q^{L_0})z^{m-|a|}w^{-m-|b|}\\
&\phantom{\phi_M Y(u,z)}=q^{-|b|}\phi_M Y(b,wq^{-1})Y(a,z)q^{L_0}.
\end{split}
\end{equation}
%\qed
\end{proof}

\begin{definition}
Let $V$ be a vertex operator algebra of central charge $c_V$
satisfying Zhu's finiteness condition.
Let $M$ be a $V$-module, $P$ a projective commutant of $M$
and $\phi$ a symmetric linear function on $P$.
We define the \textit{pseudo-trace function} 
$S_M^{P,\phi}$ on $M$ by
\begin{equation}
S_M^{P,\phi}(a,\tau)=\sum_{n=0}^\infty \phi_{M_{(r+n)}}(J_0(a)q^{L_0-c_V/24})\qquad(q=e^{2\pi i\tau},\,\tau\in\mathscr{H})
\end{equation}
for any $a\in V$.
\end{definition}

By using Proposition \ref{prop-cyclic} and the discussions given in 
\cite{Z1}, we have:

\begin{proposition}
[{\cite[Proposition 4.3.5, Proposition 4.3.6, Lemma 4.4.3]{Z1}}]
\label{prop-gtf-one-point}
A pseudo-trace function $S_M^{P,\phi}$ satisfies
\begin{align}
&S_M^{P,\phi}(a_{[0]}b,\tau)=0,\label{eq-0406-1}\\
&S_M^{P,\phi}(a_{[-2]}b,\tau)+\sum_{k=2}^\infty(2k-1)G_{2k}(\tau)
S_M^{P,\phi}(a_{[2k-2]}b,\tau)=0,\\
&S_M^{P,\phi}(L_{[-2]}a,\tau)=(2\pi i)^2q\frac{d}{dq}S_M^{P,\phi}(a,\tau)
+\sum_{k=1}^\infty G_{2k}(\tau)S_M^{P,\phi}(L_{[2k-2]}a)\label{eq-0406-2}
\end{align}
for all $a, b\in V$.
\end{proposition}

It is well known that any function subject to \eqref{eq-0406-1}--\eqref{eq-0406-2}
is a formal solution of  the ordinary differential equation 
with a regular singularity at $q=0$:
\begin{equation}
\Bigl(q\frac{d}{dq}\Bigr)^sS_M^{P,\phi}(a,\tau)+\sum_{i=0}^{s-1}h_i(q)
\Bigl(q\frac{d}{dq}\Bigr)^i
S_M^{P,\phi}(a,\tau)=0,
\end{equation}
where $h_i\in\C[G_2,G_4,G_6]$ (see \cite{Z1}, \cite{DLM3}).
Since quasi-modular forms $h_i(q)$ converge on the domain $|q|<1$, 
$S_M^{P,\phi}(a,\tau)$ converges on the same domain.
Therefore the pseudo-trace function $S_M^{P,\phi}(a,\tau)$ is holomorphic on $\mathscr{H}$.

Proposition \ref{prop-gtf-one-point} 
together with the above discussions shows:
\begin{theorem}\label{thm-gtf}
Let $V$ be a vertex operator algebra satisfying 
Zhu's finiteness condition
and let $M=\bigoplus_{n=0}^\infty M_{(r+n)}$ be a $V$-module.
For a projective commutant $P$ of $M$ and a symmetric linear function $\phi$
on $P$, the pseudo-trace function $S_M^{P,\phi}$ is a one-point function.
\end{theorem}

\begin{remark}
Let $V$ be a vertex operator algebra satisfying Zhu's finiteness condition.
Let $M$ be a simple $V$-module and set $P=\End_V(M)\cong\C$.
Then the pseudo-trace function $S_M^{P,\,\id_P}$ is noting but 
the usual trace function on $M$ which gives rise to the character of $M$
when it is evaluated at the vacuum vector.
\end{remark}

Let $P$ be a projective commutant of a $V$-module 
$M=\bigoplus_{n=0}^\infty M_{(r+n)}$.
Recall that the space of symmetric linear functions $\slf^P$ is a $Z(P)$-module. Then, by the definition of pseudo-trace functions, we have:

\begin{proposition}\label{prop-centermodule}
Let $V$ be a vertex operator algebra satisfying Zhu's finiteness condition 
and let $P$ be a projective commutant of a $V$-module 
$M=\bigoplus_{n=0}^\infty M_{(r+n)}$.
Then the vector space linearly spanned by 
$\{S_M^{P,\phi}\,|\,\phi\in\slf^P\}$ is a $Z(P)$-module by the action defined by
\begin{equation}
c\cdot S_M^{P,\phi}=S_M^{P,c\cdot\phi}.
\end{equation}
\end{proposition}

\section{The vertex operator algebra $\mathscr{F}^+$}\label{sect-5}
In this section, we recall from \cite{A} the definition of the $\Z_2$-orbifold model
$\mathscr{F}^+$ of the symplectic fermionic vertex operator superalgebra
and the classification of its simple modules.
We also review the construction of certain indecomposable modules
and the structure of  Zhu's algebra for $d=1$, which are given in \cite{A}.
\subsection{The vertex operator algebra $\mathscr{F}^+$ and its simple $\mathscr{F}^+$-modules}\label{sect-5-1}
Let $\frakh$ be a $2d$-dimensional $\C$-vector space
with a non-degenerate skew symmetric bilinear form 
$\langle\;,\,\rangle:\frakh\times\frakh\to\C$.
Then there exists a basis $\{\phi^i,\psi^i\}_{i=1}^d$ of $\frakh$ such that 
\begin{equation}
\langle \phi^i,\phi^j\rangle=\langle \psi^i,\psi^j\rangle=0,\quad
\langle \phi^i,\psi^j\rangle=-\langle \psi^j,\phi^i\rangle=-\delta_{ij}.
\end{equation}
We denote by $\hat{\frakh}$ the affinization of $\frakh$, that is,
\begin{equation}
\hat{\frakh}=\frakh\otimes\C[t,t^{-1}]\oplus\C K.
\end{equation}
The vector space $\hat{\frakh}$ becomes a superspace by letting 
$\C K$ be an even part and
$\frakh\otimes\C[t,t^{-1}]$ an odd part.
The commutation relations on $\hat{\frakh}$ are given by
\begin{equation}\label{4.1}
[a\otimes t^m,b\otimes t^n]_+=m\langle a,b\rangle\delta_{m+n,0}K,\quad
[K,\hat{\frakh}]=0
\end{equation}
for all $a, b\in\frakh$ and $m, n\in\Z$,
where $[\,,\,]_+$ denotes the anti-commutator.

We set $\mathscr{A}=U(\hat{\frakh})/(K-1)$ where $U(\hat{\frakh})$
is the universal enveloping algebra of the Lie superalgebra $\hat{\frakh}$ 
and $(K-1)$ is the two-sided ideal of $U(\hat{\frakh})$ generated 
by $K-1$. We denote $h\otimes t^n$ by $h_{n}$ for any integer $n$.
The $\Z_2$-grading on $\hat{\frakh}$ induces the $\Z_2$-grading on
$\mathscr{A}$, that is, $\mathscr{A}$ decomposes into a direct sum 
of the even part $\mathscr{A}_{\bar{0}}$ 
and the odd part $\mathscr{A}_{\bar{1}}$ of $\mathscr{A}$.
More precisely, the even part $\mathscr{A}_{\bar{0}}$ is linearly spanned by
the set
$\left\{h^1_{-n_1}h^2_{-n_2}\dotsb h^{2r}_{-n_{2r}}1\,|\,r\in\Z_{\geq0}\,,h^j\in\frakh,\,n_j\in\Z\right\}$
and the odd part $\mathscr{A}_{\bar{1}}$ is linearly spanned by the set
$\left\{h^1_{-n_1}h^2_{-n_2}\dotsb h^{2r+1}_{-n_{2r+1}}1\,|\,r\in\Z_{\geq0}, h^j\in\frakh, n_j\in\Z\right\}$.

Let $\mathscr{A}_{\geq 0}$ be the left ideal of $\mathscr{A}$ generated 
by $h_{n}1$ for all $h\in\frakh$ and non-negative integers $n$.
We set $\mathscr{F}=\mathscr{A}/\mathscr{A}_{\geq 0}$
and denote the image of the unity $1\in \mathscr{A}$ by $\vac$.
Then $\mathscr{F}$ is a left $\mathscr{A}$-module and is $\Z_2$-graded since
$\mathscr{F}=(\mathscr{A}_{\geq 0}\cap\mathscr{A}_{\bar{0}})\oplus(\mathscr{A}_{\geq 0}\cap\mathscr{A}_{\bar{1}})$.

For any $h\in\frakh$, we define a field $h(z)=\sum_{n\in\Z}h_nz^{-n-1}$
on $\mathscr{F}$. Then the first of commutation relations \eqref{4.1}
is equivalent to the operator product expansion 
\begin{equation}
h^1(z)h^2(w)\sim\frac{\langle h^1,h^2\rangle}{(z-w)^2}\quad(a,\,b\in\frakh),
\end{equation}
in particular, any two fields $h^1(z)$ and $h^2(z)$ are mutually local.

For any field $h(z)\,(h\in\frakh)$, we set
\begin{equation}
h(z)_-=\sum_{n<0}h_nz^{-n-1},\quad h(z)_+=\sum_{n\geq0}h_nz^{-n-1}.
\end{equation}
We define the normally ordered product $\NO h^1(z)h^2(z) \NO$ of two fields
$h^1(z)$ and $h^2(z)$ for $h^1, h^2\in\frakh$ by
$\NO h^1(z)h^2(z) \NO=h^1(z)_-h^2(z)-h^2(z)h^1(z)_+$.
Then we extend it to
$\NO h^1(z)h^2(z)\dotsb h^n(z)\NO$, recursively, that is,
\begin{equation}
\NO h^1(z)h^2(z)\dotsb h^n(z)\NO=
\NO h^1(z)\NO h^2(z)\dotsb h^n(z)\NO\NO.
\end{equation}
Then we obtain a linear map 
$Y:\mathscr{F}\mapsto\End_\C(\mathscr{F})[[z,z^{-1}]]$
which is defined by 
\begin{equation}
Y(h^1_{-n_1-1}\dotsb h^r_{-n_r-1}\vac,z)
=\NO\partial^{(n_1)}h^1(z)\dotsb \partial^{(n_r)}h^r(z)\NO,
\end{equation}
where $\partial^{(n)}=\frac{1}{n!}\frac{d^n}{dz^n}$ and $n_i\,(1\le i\le r)$ are non-negative integers.
Then we set
$Y(a,z)=\sum_{n\in\Z}a_{(n)}z^{-n-1}$ for any $a\in\mathscr{F}$.

By the definition of linear maps
$a_{(n)}\,(a\in\mathscr{F}, n\in\Z)$, it follows that $(h_{-1}\vac)_{(n)}=h_{n}$
for any $h\in\frakh$ and $n\in\Z$.
Therefore the superspace $\frakh$ is identified with a super subspace
of $\mathscr{F}$
by the injective map $\frakh\rightarrow\mathscr{F}\,
(h\mapsto h_{-1}\vac)$.

We set $\omega=\sum_{i=1}^d\phi^i_{-1}\psi^i$ and denote 
the corresponding vertex operator by 
$L(z)=\sum_{n\in\Z}L_nz^{-n-2}$.
Since 
\begin{equation}
L(z)L(w)
\sim-\frac{d}{(z-w)^4}+\frac{2L(w)}{(z-w)^2}
+\frac{\partial L(w)}{z-w},
\end{equation}
the operators $L_{n}=\omega_{(n+1)}$ with $n\in\Z$ as well as the identity map $\id_{\mathscr{F}}$
give rise to a module structure of the Virasoro algebra
of central charge $-2d$ on $\mathscr{F}$.
Moreover we have
\begin{equation}
\begin{split}
L(z)h(w)&\sim\sum_{i=1}^d\NO\phi^i(z)\psi^i(z)\NO h(w)\\
&\sim\sum_{i=1}^d\left(-\langle\phi^i(z)h(w)\rangle\psi^i(z)
+\langle\psi^i(z)h(w)\rangle\phi^i(z)\right)\\
&\sim\sum_{i=1}^d\left(-\frac{\langle\phi^i,h\rangle\psi^i(z)}{(z-w)^2}+
\frac{\langle\psi^i,h\rangle\phi^i(z)}{(z-w)^2}\right)\\
&\sim\frac{h(z)}{(z-w)^2}\\
&\sim\frac{h(w)}{(z-w)^2}+\frac{\partial h(w)}{z-w}.
\end{split}
\end{equation}
Then
\begin{equation}\label{weight}
[L_m,h_{n}]=-nh_{m+n}
\end{equation}
for all $m, n\in\Z$. In particular, $[L_{-1},h_{n}]=-nh_{n-1}$ 
for all integers $n$, which is equivalent to
$[L_{-1},h(z)]=\frac{d}{dz} h(z)$. 
By the definition of the field $L(z)$, we have 
\begin{equation}
L_0h=h,\quad L_nh=0
\end{equation}
for any $h\in\frakh$ and any positive integer $n$.

By using the existence theorem (cf. \cite[Theorem 4.5]{K}), we have:

\begin{theorem}{\rm(\cite[Theorem 3.1]{A})}\label{4-1-1}
The quadruple $(\mathscr{F},Y,\vac,\omega)$ is 
a simple vertex operator superalgebra of central charge $-2d$ 
with the vacuum vector $\vac$ and the Virasoro vector $\omega$.
The gradation of  $\mathscr{F}$ is given by 
$\mathscr{F}=\bigoplus_{n=0}^{\infty}\mathscr{F}_{n}$
where $\mathscr{F}_{n}$ is linearly panned by
\begin{equation}
h^{1}_{-n_1}h^2_{-n_2}\dotsb h^{r}_{-n_r}\vac\qquad(n_1+n_2+\dotsb+n_r=n)
\end{equation}
with $r\in\Z_{\ge0}$, $h^j\in\frakh$ and $n_j\in\Z_{>0}$ 
\end{theorem}

It is obvious that  the even part $\mathscr{F}^+$ of 
$\mathscr{F}$ is a vertex operator algebra
and the odd part $\mathscr{F}^-$ is a simple module 
for $\mathscr{F}^+$:
  
\begin{proposition}[{\cite[Proposition 3.2]{A}}]\label{4-1-2}
The even part $\mathscr{F}^+$ of $\mathscr{F}$ 
is a simple vertex operator algebra of central charge $-2d$ 
with the vacuum vector $\vac$ and the Virasoro vector $\omega$.
The odd part $\mathscr{F}^-$ of $\mathscr{F}$ 
is a simple $\mathscr{F}^+$-module.
\end{proposition}

Recall that $\mathscr{F}^+=\bigoplus_{n=0}^{\infty}\mathscr{F}_{n}^+$
where $\mathscr{F}^+_{n}$ is linearly spanned by elements
\begin{equation}
h^{1}_{-n_1}h^2_{-n_2}\dotsb h^{2r}_{-n_{2r}}\vac\qquad(n_1+n_2+\dotsb+n_{2r}=n)
\end{equation}
with $r\in\Z_{\ge0}$, $h^j\in\frakh$ and $n_j\in\Z_{>0}$.
In particular, we have
$\mathscr{F}^+_0=\C\vac$ and $\mathscr{F}^+_1=0$.
We also have $\mathscr{F}^-=\bigoplus_{n=1}^{\infty}\mathscr{F}^-_{n}$
where $\mathscr{F}^-_{n}$ is linearly spanned by elements
\begin{equation}
h^{1}_{-n_1}h^2_{-n_2}\dotsb h^{2r+1}_{-n_{2r+1}}\vac\qquad(n_1+n_2+\dotsb+n_{2r+1}=n)
\end{equation}
with $r\in\Z_{\ge0}$, $h^j\in\frakh$ and $n_j\in\Z_{>0}$.

\begin{remark}
The vertex operator algebra $\mathscr{F}^+$ for $d=1$ is isomorphic to 
the triplet $\mathscr{W}$-algebra $\mathscr{W}(2)$ 
(see \cite[Remark 3.9]{A} and references therein).
\end{remark}

One of the important features of the vertex operator algebra $\mathscr{F}^+$ is:
\begin{theorem}[{\cite[Theorem 3.10]{A}}]\label{4-1-3}
The vertex operator algebra $\mathscr{F}^+$ satisfies Zhu's finiteness condition. 
\end{theorem}

As we have already mentioned, there are two inequivalent simple 
$\mathscr{F}^+$-modules $\mathscr{F}^+$ and $\mathscr{F}^-$.
In \cite{A},  two  more simple $\mathscr{F}^+$-modules
$\mathscr{F}_{t}^{+}$ and $\mathscr{F}_{t}^{\,-}$ are found
as the even and the odd parts of the simple twisted $\mathscr{F}$-module
$\mathscr{F}_t$.
The conformal weights of $\mathscr{F}_{t}^{+}$ and $\mathscr{F}_{t}^{\,-}$ are $-d/8$ and $(-d+4)/8$, respectively.

\begin{theorem}[{\cite[Theorem 4.2]{A}}]\label{4-1-4}
The complete list of inequivalent simple $\mathscr{F}^+$-modules is
$\{\mathscr{F}^{\pm},\mathscr{F}_{t}^{\pm}\}$. 
\end{theorem}

\begin{corollary}\label{4-1-5}
The set of conformal weights is 
$\Lambda=\left\{0,\,1,\,-d/8,\,(-d+4)/8\right\}$.
\end{corollary}

\subsection{Indecomposable modules for the vertex operator algebra $\mathscr{F}^+$}\label{sect-5-2}

Let $\mathscr{A}_+$ be the left ideal of $\mathscr{A}$ 
which is generated by $h_{n}1$ with $h\in\frakh$ and $n\in\Z_{>0}$.
Then $\mathscr{F}_+=\mathscr{A}/\mathscr{A}_+$ is 
a left $\mathscr{A}$-module. We denote the image of 
$1$ in $\mathscr{A}$ by $\vac_+$. The left $\mathscr{A}$-module 
$\mathscr{F}_+$ is isomorphic to
the exterior algebra $\bigwedge(\frakh\otimes\C[t^{-1}])$ as vector spaces
(see \cite{A}).

\begin{proposition}[{\cite[Section 5]{A}}]
The left $\mathscr{A}$-module $\mathscr{F}_+$ is an $\mathscr{F}$-module
whose  module structure is given by 
\begin{equation}\label{eq-vertex-operator-indec}
Y(h^1_{-n_1-1}\dotsb h^r_{-n_r-1}\vac,z)
=\NO\partial^{(n_1)}h^1(z)\dotsb \partial^{(n_r)}h^r(z)\NO.
\end{equation}
\end{proposition}

It follows from the construction of $\mathscr{F}_+$ 
that the $\mathscr{F}^+$-module $\mathscr{F}_+$ is linearly spanned 
by elements 
\begin{equation}\label{eq-indecomp-spanning}
h^1_{-n_1}h_{-n_2}^2\dotsb h^r_{-n_r}\vac_+\quad(r\in\Z_{\ge0}, h^j\in\frakh, n_j\in\Z_{\geq0}).
\end{equation}
Since
\begin{equation}\label{eq-gen1025}
\begin{split}
L_0\vac_+&=\sum_{i=1}^d(\phi^i_{-1}\psi^i_{-1}1)_{(1)}\vac_+\\
&=\sum_{i=1}^d\Bigl(\sum_{n<0}\phi^i_{n}\psi^i_{-n}-\sum_{n\geq 0}\psi^i_{-n}\phi^i_{n}\Bigr)\vac_+\\
&=\sum_{i=1}^d\phi^i_{0}\psi^i_{0}\vac_+
\end{split}
\end{equation}
and $(h_{0})^2=0$ for any $h\in\frakh$, we have 
\begin{equation}\label{eq-vacuumlike}
L_0^{d+1}\vac_+=0.
\end{equation}
Then, by \eqref{weight}, we see that $\mathscr{F}_+$
decomposes into a direct sum of finite-dimensional generalized eigenspaces for $L_0$,
that is, $\mathscr{F}_+=\bigoplus_{n=0}^\infty \mathscr{F}_{+(n)}$ where
where $\mathscr{F}_{+(n)}$
is linearly spanned by elements
\begin{equation}\label{eq-spanning}
h^1_{-n_1}h^2_{-n_2}\dotsb h^r_{-n_r}
\bigl(\prod_{i=1}^d(\phi_0^i)^{m_i}(\psi_0^i)^{n_i}\bigr)\vac_+\qquad(m_i,n_i=0,1)
\end{equation}
with $r\in\Z_{\ge0}$, $h^j\in\frakh$, $n_j\in\Z_{>0}$
and $\sum_{j=1}^{r}n_j=n$
and $(L_0-n)^{d+1}\mathscr{F}_{+(n)}=0$.

Since the $\Z_2$-grading of $\hat{\frakh}$ also induces 
a $\Z_2$-grading on $\mathscr{A}_+$,
it follows that $\mathscr{F}_+$ is a superspace with
the even part $\mathscr{F}_+^+$
and the odd part $\mathscr{F}_+^-$.
Note that $\mathscr{F}_+^+$ is linearly spanned by elements of the form \eqref{eq-indecomp-spanning} of even length
and $\mathscr{F}_+^-$ is linearly spanned by elements of the form \eqref{eq-indecomp-spanning} of odd length.

\begin{proposition}[{\cite[Corollary 5.2]{A}}]\label{prop-0406}
The $\mathscr{F}^{+}$-modules $\mathscr{F}_+^{\pm}$ are reducible and indecomposable. 
\end{proposition}

\begin{remark}
The vertex operator algebra $\mathscr{F}^+$ satisfies Zhu's finiteness condition
but it is not rational by Proposition \ref{prop-0406}.
\end{remark}
\subsection{Zhu's algebra $A_0(\mathscr{F}^+)$ for $d=1$}\label{sect-5-3}
Let us denote by $\mathscr{T}^+$ the vertex operator algebra $\mathscr{F}^+$ for $d=1$.
The minimal polynomial of $[\omega]$ of 
$A_0(\mathscr{T}^+)$ is given as follows:

\begin{proposition}[{\cite[Proposition 4.4]{A}}]
\label{prop-zhu-relation}
The minimal polynomial of the central element $[\omega]$ of $A_0(\mathscr{T}^+)$ is
\begin{equation}
[\omega]^2*([\omega]-1)*(8[\omega]+1)*(8[\omega]-3)=0.
\end{equation}
\end{proposition}

The proposition yields
\begin{equation}
A(\mathscr{T}^+)
=A_0\oplus A_1\oplus A_{-\frac{1}{8}}\oplus A_{\frac{3}{8}}
\end{equation}
which is the decomposition into a direct sum of two-sided ideals,
where $A_\lambda$ is the generalized eigenspace of $[\omega]$
with eigenvalue $\lambda$.

The structure of Zhu's algebra 
$A_0(\mathscr{T}^+)$
is described as follows.

\begin{proposition}[{\cite[Proposition 4.6]{A}, \cite[Theorem 4.6]{NT2}}, \cite{AM2}]\label{prop-shape}
\begin{enumerate}
\item[{\rm (1)}]
For $\lambda=1$ and $\lambda=3/8$, the two-sided ideal $A_\lambda$
is isomorphic to the $2\times 2$ matrix algebra $M_2(\C)$.
\item[{\rm (2)}]
The two-sided ideal $A_{-1/8}$ is one-dimensional.
\item[{\rm(3)}]
The two-sided ideal $A_0$ is isomorphic to the algebra
\begin{equation}
\left\{
\begin{pmatrix}
a&b\\
0&a
\end{pmatrix}\,\Big|\,a,b\in\C\right\}.
\end{equation}
\end{enumerate}
\end{proposition}

Since  the vector space of symmetric linear functions on a matrix algebra
is one-dimensional and $A_0$ is a commutative algebra,
we have the following:
\begin{corollary}\label{cor-slf-dim-d1}
The vector space $\slf^{\mathscr{T}^+}$ is $5$-dimensional.
\end{corollary}

\section{Pseudo-trace functions on $\mathscr{F}_+$}\label{sect-6}
In this section, we determine the endomorphism ring of $\mathscr{F}_+$
and prove that the endomorphism ring is itself a projective commutant of $\mathscr{F}_+$.
We also determine the vector space of symmetric linear functions on the endomorphism ring.
Then we introduce a $P$-coordinate system of $\mathscr{F}_+$
and construct pseudo-trace functions on $\mathscr{F}_+$.

\subsection{The endomorphism ring of $\mathscr{F}_+$}\label{sect-6-1}
In order to construct pseudo-trace functions on $\mathscr{F}_+$,
we need a projective commutant of $\End_{\mathscr{F}^+}(\mathscr{F}_+)$
and symmetric linear functions on it.

Recall that $\mathscr{F}^+$ is linearly spanned by ``even length elements''
\begin{equation}
h^{1}_{-n_1}h^2_{-n_2}\dotsb h^{2r}_{-n_{2r}}\vac\qquad(n_1+n_2+\dotsb+n_{2r}=n)
\end{equation}
with $r\in\Z_{\ge0}$, $h^j\in\frakh$ and $n_j\in\Z_{>0}$.
By \eqref{eq-vertex-operator-indec},
the operators $a_{(n)}\,(a\in\mathscr{F}^+, n\in\Z)$ on $\mathscr{F}_+$ are written
as linear combinations of operators $h^{1}_{n_1}h^2_{n_2}\dotsb h^{2r}_{n_{2r}}$
with $h^j\in\frakh$ and $n_j\in\Z$.
Hence, by \eqref{4.1}, it follows that operators $\phi^i_{0}$ and $\psi_{0}^i\,(1\leq i\leq d)$ 
on $\mathscr{F}_+$ commute with $a_{(n)}$ for any $a\in\mathscr{F}^+$
and integers $n$.
Therefore it follows that 
$\phi^i_{0},\,\psi^i_{0}\in\End_{\mathscr{F}^+}(\mathscr{F}_+)$. 
Let $\theta$ be the automorphism
on $\mathscr{F}_+=\mathscr{F}_+^+\oplus\mathscr{F}_+^-$ defined by 
$a+b\mapsto a-b$ where $a\in\mathscr{F}_+^+$ 
and $b\in\mathscr{F}_+^-$, respectively.
Let $P$ be the subalgebra of $\End_{\mathscr{F}^+}(\mathscr{F}_+)$ which is generated by 
$e_i=\phi^i_{0}, f_i=\psi^i_{0}$ and $K=\theta$.
Noting that $h^1_{0}h^2_{0}+h^2_{0}h^1_{0}=0$ 
and $\theta h^1_{0}+h^1_{0}\theta=0$ for any $h^1, h^2\in\frakh$,
we have
\begin{equation}\label{eq-fund-relation}
\begin{split}
&K^2=1,
\quad e_ie_j=-e_je_i,
\quad f_if_j=-f_jf_i,\\
&Ke_i=-e_iK,
\quad Kf_i=-f_iK,\quad e_if_j=-f_je_i
\end{split}
\end{equation}
for $1\le i, j\le d$.
Note that $e_i^2=f_i^2=0$ for all $1\le i \le d$.

Now we can describe a basis of the algebra $P$.

\begin{proposition}
The set 
\begin{equation}
\Omega
=\left\{\Bigl(\prod_{i=1}^de_i^{m_i}f_i^{n_i}\Bigr)K^\ell\,|\,m_i,\,n_i,\,\ell=0,\,1
\right\}
\end{equation}
is a basis of $P$ where $\Bigl(\prod_{i=1}^de_i^{m_i}f_i^{n_i}\Bigr)K^\ell
=(e_1^{m_1}f_1^{n_1}e_2^{m_2}f_2^{n_2}\dotsb)K^\ell$.
\end{proposition}
\begin{proof}
The definition of $P$ and \eqref{eq-fund-relation} imply that 
$P$ is linearly spanned by $\Omega$.

Suppose that $\sum_{c\in\Omega}a_c c=0$
where $a_c$ are complex numbers.
Let us set the subset $\overset{\circ}{\Omega}=\left\{\prod_{i=1}^de_i^{m_i}f_i^{n_i}\,|\,m_i, n_i=0, 1\right\}$
of $\Omega$.
Then we have
\begin{equation}\label{eq-basis-proof}
0=\sum_{c\in\Omega}a_c c=\sum_{c\in\overset{\circ}{\Omega}}(a_c c+a_{cK} cK).
\end{equation}
Applying \eqref{eq-basis-proof} to $\vac_+$,
we have 
\begin{equation}\label{eq-basis-proof-2}
a_c+a_{cK}=0
\end{equation}
for all $c\in\overset{\circ}{\Omega}$
since the set $\bigl\{\prod_{i=1}^d(\phi^i_0)^{m_i}(\psi^i_0)^{n_i}\vac_+\,|\,m_i,\,n_i=0,1\bigr\}$ 
is a basis of $\mathscr{F}_{+,(0)}$.
On the other hand, applying \eqref{eq-basis-proof} to $h_{-1}\vac_+$ for 
$h\in\{\phi^i,\,\psi^i\,|\,1\leq i\leq d\}$,
we have
\begin{equation}
a_c-a_{cK}=0
\end{equation}
for any $c\in\overset{\circ}{\Omega}$
since the set
\begin{equation}
\Bigl\{\prod_{i=1}^d(\phi_0^i)^{m_i}(\psi_0^i)^{n_i}h_{-1}\vac_+\,|\,
h=\phi^i,\psi^i, m_i, n_i=0,1\Bigr\}
\end{equation}
forms a basis of $\mathscr{F}_{+,(1)}$.
Then we have $a_c=a_{cK}=0$.
%\qed
\end{proof}

We can find a symmetric linear function on $P$ which induces a non-degenerate
bilinear form on $P$.

\begin{proposition}\label{proposition5.1.1}
Let us define the linear function $\varphi:P\rightarrow\C$ by
\begin{equation}
\varphi\left(\left(\prod_{i=1}^de_i^{m_i}f_i^{n_i}\right)K^\ell\right)
=\prod_{i=1}^d\delta_{{m_i},1}\delta_{{n_i},1}\delta_{{\ell},1}.
\end{equation}
Then $P$ is a symmetric algebra with $\varphi$, that is,
$\varphi$ is symmetric and induces a non-degenerate bilinear form on $P$.
\end{proposition}
\begin{proof}
If $\varphi(ab)=0$ for $a,b\in\Omega$, then,
by \eqref{eq-fund-relation}, we have $\varphi(ab)=\pm\varphi(ba)=0$, which yields $\varphi(ab)=\varphi(ba)$.
Therefore we can assume that $\varphi(ab)\neq0$.

Suppose that $a=a_1a_2\dotsb a_\ell\in\Omega$ 
where $a_j\in\{e_i, f_i, K\,|\,1\leq i\leq d\}$.
Then there exists a unique element $b\in\Omega$ such that
$ab=\pm\left(\prod_{i=1}^de_if_i\right)K$.
The element $b$ is expressed as $b=b_{\ell+1}b_{\ell+2}\cdots b_{2d+1}$
with  $b_j\in\{e_i, f_i, K\,|\,1\leq i\leq d\}\backslash\{a_1,\dotsc,a_\ell\}$.
Then we have
\begin{equation}
ab=(-1)^{\ell(2d+1-\ell)}ba=ba,
\end{equation}
which shows that $\varphi$ is symmetric.

The matrix $A=(A_{ab})$ where $A_{ab}=\varphi(ab)$ 
with $a, b\in\Omega$ is invertible
since for any $a\in\Omega$, there exists a unique $b\in\Omega$
such that $A_{ab}=\pm \delta_{ab}$.
Therefore we see that the bilinear form $P\times P\to \C$ defined by $(a, b)\mapsto\varphi(ab)$
is non-degenerate.
%\qed
\end{proof}

By the relations \eqref{eq-fund-relation}, we can determine the center of the algebra $P$.

\begin{proposition}
The center $Z(P)$ of $P$ has the basis $Z(\Omega)\subset\Omega$ 
which consists of monomials of $e_i$ and $f_i\,(1\leq i\leq d)$ of even length,
and $\left(\prod_{i=1}^de_if_i\right)K$.
In particular, $\dim_\C Z(P)=2^{2d-1}+1$.
\end{proposition}

The map $Z(P)\rightarrow\slf^P$
defined by $c\rightarrow c\cdot\varphi$ 
where $c\cdot\varphi(x)=\varphi(cx)$ for any $x\in P$ is an isomorphism
since $P$ is a symmetric algebra with $\varphi$ (see \cite[Lemma 2.5]{Brou}). 

\begin{corollary}\label{}
We have
$\slf^P=\{c\cdot\varphi\,|\, c\in Z(P)\}$, in particular, 
$\dim_\C\slf^P=\dim_\C Z(P)=2^{2d-1}+1$.
\end{corollary}

\subsection{$P$-coordinate systems of $\mathscr{F_+}$}\label{sect-6-2}
In order to construct  pseudo-trace functions on $\mathscr{F}_+$
from symmetric linear functions on $P$, we will introduce a $P$-coordinate
system of the left $P$-module $\mathscr{F}_+$. 

Set $T^\pm=(1\pm K)/2$ and let $Q^\pm$ be $P$-submodules generated by 
$T^\pm$, respectively.
We can show that $T^\pm$ are primitive idempotents of $P$;
suppose that $M$ is a simple $P$-submodule of $Q^\pm$.
For any non-zero element $x$ of $M$, there exists an element $p$ of $P$
such that $px=\bigl(\prod_{i=1}^de_if_i\bigr)T^\pm$. 
So the submodule $M$ contains a simple $P$-submodule $\C\bigl(\prod_{i=1}^de_if_i\bigr)T^\pm$.
Since $Q^\pm$ have unique simple submodules, it follows that 
$Q^\pm$ are indecomposable.
Therefore $T^\pm$ are primitive idempotents. 
Thus $\{Q^\pm\}$ gives rise to the complete list
of indecomposable projective $P$-modules.

By \eqref{eq-spanning}, we see that any homogeneous subspace 
$\mathscr{F}_{+(n)}$
of $\mathscr{F}_+$ is a direct sum of the vector subspaces with a basis
\begin{equation}\label{eq-basis}
\Bigl\{u_{\vec{n}}^r\prod_{i=1}^d(\phi_{0}^i)^{k_i}(\psi_{0}^i)^{\ell_i}\vac_+\,|\,k_i, \ell_i=0, 1\Bigr\}
\end{equation}
for fixed $u_{\vec{n}}^r=u^1_{-n_1}u^2_{-n_2}\dotsb u^r_{-n_r}$ where 
$n_1\geq\dotsc\geq n_r>0,\,\sum_{j=1}^rn_j=n,\,u^j\in\{\phi,\,\psi\}$
and $n_j\ne n_{j+1}$ when $u^j=u^{j+1}$. 
The vector subspace linearly spanned  by the set
\eqref{eq-basis}  is  a $P$-module.
If $r$ is even, this is isomorphic to $Q^+$
as left $P$-modules by the $P$-homomorphism induced by the correspondence
\begin{equation}\label{eq-coordinate+}
u^r_{\vec{n}}\vac_+\mapsto T^+.
\end{equation}
On the other hand, if $r$ is odd, this is isomorphic to $Q^-$ 
as left $P$-modules by the $P$-homomorphism induced by
\begin{equation}\label{eq-coordinate-}
u^r_{\vec{n}}\vac_+\mapsto T^-.
\end{equation}
Since the vector subspace of $\mathscr{F}_+$ which is linearly spanned by
elements of the form 
\begin{equation}
h^1_{-n_1}h^2_{-n_2}\dotsb h^{r}_{-n_{r}}\prod_{i=1}^d\phi_{0}^i\psi_{0}^i\vac_+
\quad(r\in\Z_{\ge0}, h^j\in\frakh,\,n_j\in\Z_{>0})
\end{equation}
is isomorphic to $\mathscr{F}=\mathscr{F}^+\oplus\mathscr{F}^-$ 
as $\mathscr{F}^+$-modules, it follows that
\begin{equation}
\mathscr{F}_{+(n)}\cong
(Q^+)^{\oplus\dim_\C\mathscr{F}^+_{n}}
\oplus(Q^-)^{\oplus\dim_\C\mathscr{F}^{\,-}_{n}}.
\end{equation}
\begin{proposition}\label{prop-20111022-1}
The set
\begin{equation}
\{v_{n,i}^+, v_{n,j}^-, \alpha_{n,i}^+, \alpha^-_{n,j}\,|\,1\leq i\leq\dim_\C\mathscr{F}^+_{n},
1\leq j\leq\dim_\C\mathscr{F}^-_{n}\}
\end{equation}
is a $P$-coordinate system of $\mathscr{F}_{+(n)}$ where 
$\alpha_{n,i}^+$ is defined by \eqref{eq-coordinate+},
$\alpha_{n,j}^-$ is defined by \eqref{eq-coordinate-}, and
$v_{n,i}^+=u^r_{\vec{n}}\vac_+$ for an even  $r$,
$v_{n,i}^-=u^r_{\vec{n}}\vac_+$ for an  odd $r$.
\end{proposition}
\begin{proof}
Let $Q_{n,i}^\pm$ be left $P$-submodules of $\mathscr{F}_{+(n)}$ generated by $v_{n,i}^\pm=u^r_{\vec{n}}\vac_+$, respectively.
Then the modules $Q^\pm_{n,i}$ are isomorphic to $Q^\pm$, respectively.
For any monomial $x$ of $e_i$ and $f_i\,(1\leq i\leq d)$, we have
\begin{equation}
\alpha_{n,i}^\pm(u^r_{\vec{n}}x\vac_+)v_{n,i}^\pm=(-1)^{|x|}xT^\pm v_{n,i}^\pm
=(-1)^{|x|}x v_{n,i}^\pm=u^r_{\vec{n}}x\vac_+,
\end{equation}
where $|x|\in\Z_2$ indicates $u^r_{\vec{n}}x\vac_+=(-1)^{|x|}xu^r_{\vec{n}}\vac_+$.
Therefore, $\{v_{n,i}^\pm, \alpha_{n,i}^\pm\}$ are
$P$-coordinate systems of $Q^\pm$, respectively.
%\qed
\end{proof}

\subsection{Pseudo-trace functions on $\mathscr{F}_+$}\label{sect-6-3}
In this subsection, we will construct pseudo-trace functions on 
$\mathscr{F}_+$ associated with symmetric linear functions on the algebra $P$. 
For any $g\in\slf^P$, we denote 
the pseudo-trace function $S^{P,g}_{\mathscr{F}_+}$
by $S^{g}$ for short.

Set $\nu=\sum_{i=1}^de_if_i\in Z(P)$.
By \eqref{weight} and \eqref{eq-gen1025}, we have $L_0-n=\nu$ as operators on $\mathscr{F}_{+(n)}$ 
for any non-negative integer $n$. 
Note that $\nu^d=d!\prod_{i=1}^de_if_i$ and $\nu^{d+1}=0$.

Recall that $P$ is a symmetric algebra with $\varphi$.
Given a symmetric linear function $c\cdot\varphi$ with $c\in Z(P)$, 
we can define the pseudo-trace function $S^{c\cdot\varphi}$.
Since $(L_0-n)^{d+1}\mathscr{F}_{+(n)}=0$ 
for any non-negative integer $n$, we have
\begin{equation}\label{pseudo-1}
\begin{aligned}
S^{c\cdot\varphi}(a,\tau)&=\sum_{j=0}^d\sum_{n=0}^\infty
\frac{1}{j!}(c\cdot\varphi)_{\mathscr{F}_{+(n)}}
\left(J_0(a)(L_0-n)^j\right)q^{n+d/12}(2\pi i\tau)^j\\
&=\sum_{j=0}^d\sum_{n=0}^\infty
\frac{1}{j!}(c\cdot\varphi)_{\mathscr{F}_{+(n)}}
\left(J_0(a)\nu^j\right)q^{n+d/12}(2\pi i\tau)^j.
\end{aligned}
\end{equation}
for all $a\in\mathscr{F}^+$.
Each term which appears in the right-hand side of \eqref{pseudo-1} can
be rewritten by using the coordinate system given in Proposition \ref{prop-20111022-1}:
\begin{multline}\label{pseudo-2}
(c\cdot\varphi)_{\mathscr{F}_{+(n)}}\left(J_0(a)\nu^j\right)\\
=\varphi\Bigl(\,\sum_{i=1}^{\dim_\C\mathscr{F}^+_n}
\alpha_{n,i}^+(J_0(a)c\nu^j v_{n,i}^+)
+\sum_{i=1}^{\dim_\C\mathscr{F}^-_{n}}
\alpha_{n,i}^-(J_0(a)c\nu^j v_{n,i}^-)\Bigr).
\end{multline}

Suppose that $c\in Z(\Omega)$  is of length $2k$.
Then $c\nu^j=0$ for 
$j>d-k$ since $\nu^j\,(0\leq j\leq d)$ is a linear combination of
elements in $Z(\Omega)$ of length $2j$.
Hence \eqref{pseudo-2} yields
that $(c\cdot\varphi)_{\mathscr{F}_{+(n)}}\left(J_0(a)\nu^j\right)=0$
for $j>d-k$.

\begin{proposition}\label{prop-length-log}
Suppose that $c\in Z(\Omega)$ is of length $2k$.
Then the pseudo-trace function $S^{c\cdot\varphi}$
is
\begin{equation}
S^{c\cdot\varphi}(a,\tau)
=\sum_{j=0}^{d-k}\sum_{n=0}^\infty S_{jn}(a)q^{n+d/12}(2\pi i\tau)^j
\end{equation}
for all $a\in V$
where $S_{jn}\in\Hom_{\C}(\mathscr{F}^+,\C)$.
\end{proposition}
We now can state one of our main results.

\begin{theorem}\label{theorem-main}
Let $\varphi$ be the symmetric linear function on $P$ defined in 
Proposition \ref{proposition5.1.1}.
\begin{enumerate}
\item[{\rm (1)}]
For $c=\bigl(\prod_{i=1}^de_if_i\bigr)K$, we have 
\begin{equation}
S^{c\cdot\varphi}(a,\tau)
=\frac{1}{2}\tr_{\mathscr{F}^+}\bigl(J_0(a)q^{L_0+d/12}\bigr)
+\frac{1}{2}\tr_{\mathscr{F}^-}\bigl(J_0(a)q^{L_0+d/12}\bigr)
\end{equation}
for all $a\in\mathscr{F}^+$
and
\begin{equation}
S^{c\cdot\varphi}(\vac,\tau)
=\frac{1}{2}\tr_{\mathscr{F}^+}\bigl(q^{L_0+d/12}\bigr)
+\frac{1}{2}\tr_{\mathscr{F}^-}\bigl(q^{L_0+d/12}\bigr).
\end{equation}
\item[{\rm (2)}]
If $c\in Z(\Omega)$ is a monomial of $e_if_i\,(1\leq i\leq d)$ of length $2k\,(k\le d)$,
then the coefficient of $(2\pi i\tau)^{d-k}$ of  $S^{c\cdot\varphi}(a,\tau)$
is
\begin{equation}\label{main-2-1}
\frac{1}{2}\left(\tr_{\mathscr{F}^+}\bigl(J_0(a)q^{L_0+d/12}\bigr)
-\tr_{\mathscr{F}^-}\bigl(J_0(a)q^{L_0+d/12}\bigr)\right)
\end{equation}
for all $a\in\mathscr{F}^+$
and
\begin{equation}\label{main-2-2}
S^{c\cdot\varphi}(\vac,\tau)=\frac{(2\pi i\tau)^{d-k}}{2}
\left(\tr_{\mathscr{F}^+}\bigl(q^{L_0+d/12}\bigr)-
\tr_{\mathscr{F}^-}\bigl(q^{L_0+d/12}\bigr)\right).
\end{equation}
\item[{\rm (3)}]
If $c\in Z(\Omega)$ is neither a monomial of 
$e_if_i\,(1\leq i\leq d)$ nor $\prod_{i=1}^de_if_iK$,
then we have $S^{c\cdot\varphi}(\vac,\tau)=0$.
\item[{\rm (4)}]
For any $c\in Z(\Omega)$, the pseudo-trace function 
$S^{c\cdot\varphi}$ is not zero and the set 
$\{S^{c\cdot\varphi}\,|\,c\in Z(\Omega)\}$ is linearly independent.
\item[{\rm (5)}]
$\dim_\C\mathscr{C}(\mathscr{F}^+)\geq2^{2d-1}+3$.
\end{enumerate}
\end{theorem}
\begin{proof}
(1)
Note that $cv_{n,i}^+=\bigl(\prod_{k=1}^de_kf_k\bigr)v_{n,i}^+$
and $cv_{n,i}^-=-\bigl(\prod_{k=1}^de_kf_k\bigr)v_{n,i}^+$.
This implies that the vector subspaces spanned by
 $\{cv_{n,i}^\pm\,|\,n\geq0, 1\leq i\leq \dim\mathscr{F}^\pm_n\}$
are isomorphic to $\mathscr{F}^\pm$, respectively.
Now we obtain, by the definition of $\varphi$,
\begin{equation}
\varphi(\alpha_{n,i}^\pm(cv_{n,i}^\pm))=\varphi\bigl(\prod_{k=1}^de_if_i\alpha_{n,i}^\pm
(v_{n,i}^\pm)\bigr)=\varphi\bigl(\prod_{k=1}^de_if_iT^\pm\bigr)=\pm\frac{1}{2}.
\end{equation}
Therefore the sets 
$\{2\varphi\circ\alpha_{n,i}^\pm\,|\,1\leq i\leq \dim\mathscr{F}_n^\pm\}$
are dual to the bases $\{cv_{n,i}^\pm\,|\,1\le i\le\dim_{\C}\mathscr{F}^\pm_n\}$
of $\mathscr{F}^\pm_n$, respectively.
Then we have
\begin{equation}
\begin{split}
&\varphi\Bigl(\,\sum_{i=1}^{\dim\mathscr{F}^+_n}
\alpha_{n,i}^+(J_0(a)c v_{n,i}^+)
+\sum_{i=1}^{\dim\mathscr{F}^-_{n}}
\alpha_{n,i}^-(J_0(a)c v_{n,i}^-)\Bigr)\\
&\phantom{\varphi\Bigl(\,\sum_{i=1}^{\dim\mathscr{F}^+_n}\alpha_{n,i}^+(J_0(a)c v_{n,i}^+)+}
=\frac{1}{2}(\tr_{\mathscr{F}^+_n}(J_0(a))+\tr_{\mathscr{F}^-_n}(J_0(a))).
\end{split}
\end{equation}
 
Remark that  $c\nu^j=\left(\prod_{i=1}^de_if_i\right)K\nu^j=0$ 
for all positive integers $j$.
Then \eqref{pseudo-1} and \eqref{pseudo-2} show the statement (1).

(2) 
By \eqref{pseudo-1} and \eqref{pseudo-2}
the coefficient of $(2\pi i\tau)^{d-k}q^{n+d/12}$ 
of  $S^{c\cdot\varphi}(a,\tau)$ is
\begin{equation}\label{eq-5-27-co}
\frac{1}{(d-k)!}\varphi\Bigl(\;\sum_{i=1}^{\dim_\C\mathscr{F}^+_n}
\alpha_{n,i}^+(J_0(a)c\nu^{d-k} v_{n,i}^+)
+\sum_{i=1}^{\dim_\C\mathscr{F}^-_{n}}\alpha_{n,i}^-(J_0(a)c\nu^{d-k} v_{n,i}^-)\Bigr).
\end{equation}
Since $\nu^j$ is equal to
\begin{equation}\label{eq-exterior}
j!\times\{\text{the sum of all monomials of $e_if_i\,(1\leq i\leq d)$ of length $2j$}\},
\end{equation}
we have $c\nu^{d-k}=(d-k)!\prod_{i=1}^de_if_i$.
By the same argument in (1), we have
\begin{equation}
\begin{split}
&\frac{1}{(d-k)!}\varphi\biggl(\,\sum_{i=1}^{\dim_\C\mathscr{F}^+_n}
\alpha_{n,i}^+(J_0(a)c\nu^{d-k} v_{n,i}^+)
+\sum_{i=1}^{\dim_\C\mathscr{F}^-_{n}}\alpha_{n,i}^-(J_0(a)c\nu^{d-k} v_{n,i}^-)\biggr)\\
&\phantom{\frac{1}{(d-k)!}\varphi\biggl(\,\sum_{i=1}^{\dim_\C\mathscr{F}^+_n}
\alpha_{n,i}^+(J_0(a)c\nu^{d-k} v_{n,i}^+)}
=\frac{1}{2}(\tr_{\mathscr{F}^+_n}(J_0(a))-\tr_{\mathscr{F}^-_n}(J_0(a))).
\end{split}
\end{equation}
This proves the first statement of $(2)$. 
In order to prove the second statement, by virtue of Proposition \ref{prop-length-log}, it is sufficient to show
\begin{equation}\label{eq-wanted}
\varphi\biggl(\,\sum_{i=1}^{\dim_\C\mathscr{F}^+_n}
\alpha_{n,i}^+(c\nu^{j} v_{n,i}^+)
+\sum_{i=1}^{\dim_\C\mathscr{F}^-_{n}}\alpha_{n,i}^-(c\nu^{j} v_{n,i}^-)\biggr)=0
\end{equation}
for $0\leq j<d-k$.
Since $\nu^j$ is a scalar multiple of the sum of all monomials 
of $e_if_i\,(1\leq i\leq d)$
of length $2j$ and $c$ is of length $2k$,
the element  $c\nu^j$ is a scalar multiple of the sum of monomials of $e_if_i\,(1\leq i\leq d)$
whose length is smaller than $2d$ for $0\leq j<d-k$. Therefore, by the definition of $\varphi$, we have
\begin{equation}
\varphi(\alpha_{n,i}^\pm(c\nu^jv_{n,i}^\pm))=\varphi(c\nu^jT^\pm)=0,
\end{equation}
which shows \eqref{eq-wanted}.

We now give a proof of $(3)$.
Since $c\not=\bigl(\prod_{i=1}^de_if_i\bigr)K$, the length of $c\in Z(\Omega)$ is $2k$.
If $k=d$, then $c=\prod_{i=1}^de_if_i$ which contradicts to the assumption
of (3).
Then we have $k<d$.
Therefore,
by \eqref{pseudo-1}, \eqref{pseudo-2}, and Proposition \ref{prop-length-log},
in order to prove (3), 
it is sufficient to prove
\begin{equation}\label{eq-wanted-2}
\varphi(\alpha_{n,i}^\pm(c\nu^{d-k}v_{n,i}^\pm))=\varphi(c\nu^{d-k}T^\pm)=0
\end{equation} 
for all $1\leq i\leq \dim_\C\mathscr{F}^\pm_{n}$.
Since the element $c\nu^{d-k}$ is a linear combination of monomials
of length $2d$ by \eqref{eq-exterior},  the element 
$c\nu^{d-k}$ is equal 
to $\lambda\prod_{i=1}^de_if_i$ for a complex number $\lambda$.
If $\lambda$ is non-zero, then $c$ must be a monomial of $e_if_i\,(1\leq i\leq d)$.
This contradicts to the assumption of (3) and
therefore we have $\lambda=0$, that is, $c\nu^{d-k}=0$.
Hence we have shown \eqref{eq-wanted-2}. 

We will prove the statement {\rm(4)}.
By $(1)$ and $(2)$, the pseudo-trace functions 
$S^{c\cdot\varphi}$ for $c=\prod_{i=1}^de_if_i$ 
and $c=(\prod_{i=1}^de_if_i)K$ are both non-zero.
For any $c\in Z(\Omega)$ of length $2k$, there exists a unique element
$c^\prime\in Z(\Omega)$ of length $2(d-k)$ such that $cc^\prime=\pm\prod_{i=1}^de_if_i$.
By Proposition \ref{prop-centermodule} and Theorem \ref{theorem-main} (2), 
we have 
\begin{equation}\label{eq-0411}
\begin{split}
c^\prime\cdot S^{c\cdot\varphi}(a,\tau)&=\pm S^{(cc^\prime)\cdot\varphi}(a,\tau)\\
&=\pm\frac{1}{2}(\tr_{\mathscr{F}^+}(J_0(a)q^{L_0+d/12})-\tr_{\mathscr{F}^-}(J_0(a)q^{L_0+d/12}))\not=0
\end{split}
\end{equation}
for all $a\in\mathscr{F}^+$
and then $S^{c\cdot\varphi}\not=0$.

For proving the second statements of $(4)$, we suppose that
\begin{equation}\label{eq-main-4}
\sum_{c\in Z(\Omega)}a_cS^{c\cdot\varphi}=0
\end{equation}
with complex numbers $a_c$.
We proceed induction on the length of elements in $Z(\Omega)$.
For the central element $c=1$, we choose
$c^\prime=\prod_{i=1}^de_if_i$.
Then we see that $c^{\prime\prime}c^\prime=0$ for any $c^{\prime\prime}\in Z(\Omega)$ whose length
is non-zero.
Therefore we have 
\begin{equation}
c^\prime\cdot\sum_{c^{\prime\prime}\in Z(\Omega)}a_{c^{\prime\prime}}S^{c^{\prime\prime}\cdot\varphi}
=a_{c}S^{c^\prime\cdot\varphi}=0.
\end{equation}
Since $S^{c^\prime\cdot\varphi}\not=0$ we have $a_{c}=0$.

Suppose that $a_c=0$ for any $c\in Z(\Omega)$ 
whose length is strictly smaller than $2k$ with $k>0$.
Let $c$ be an element of $Z(\Omega)$ of length $2k$.
Then there exists a unique $c^\prime\in Z(\Omega)$ such that 
the length of $c^\prime$ is $2(d-k)$ and $cc^\prime=\pm\prod_{i=1}^de_if_i$. 
Note that $c^{\prime\prime}c^\prime=0$ for any $c^{\prime\prime}\in Z(\Omega)$ 
whose length is greater than or equal to $2k$ except for $c$. Hence we have 
$a_{c}=0$. This proves that 
$a_c=0$ for all $c\in Z(\Omega)$ whose length is $2k$.
Therefore, we obtain $a_c=0$ for all $c\in Z(\Omega)$ whose length is 
strictly smaller than $2d$.
We still have to show $a_c=0$ with $c=\prod_{i=1}^de_if_i$ and 
$c=\prod_{i=1}^de_if_iK$.
In these cases, by Theorem \ref{theorem-main} (1) and (2), we have
$a_c=0$.
Then we conclude that $a_c=0$ for all $c\in Z(\Omega)$.

(5) By Theorem \ref{theorem-main} (1) and (2), we obtain ordinary trace functions on simple modules
$\mathscr{F}^\pm$ as special cases of pseudo-trace 
functions. Moreover, we have two ordinary trace functions on $\mathscr{F}_t^\pm$.
If the length of $c\in Z(\Omega)$ is smaller than $2d$, then the pseudo-trace function 
$S^{c\cdot\varphi}$ has logarithmic terms.
Thus the set consisting of ordinary trace functions on $\mathscr{F}^\pm_t$
and pseudo-trace functions $S^{c\cdot\varphi}\,(c\in Z(\Omega))$
is linearly independent.
Since $|Z(\Omega)|=2^{2d-1}+1$, we see that 
$\dim_\C\mathscr{C}(\mathscr{F^+})\geq2^{2d-1}+3$ by (4).
%\qed
\end{proof}

By virtue of Theorem \ref{theorem-main} (3) and (4),
the surjective linear map
$\mathscr{C}(\mathscr{F}^+)\to \gch(\mathscr{F}^+)$
restricted to the subspace of $\mathscr{C}(\mathscr{F}^+)$
linearly spanned by $\{S^{c\cdot\varphi}\,|\,c\in Z(\Omega)\}$
is not injective. Therefore, and so is the linear map
$\mathscr{C}(\mathscr{F}^+)\to \gch(\mathscr{F}^+)$.
\begin{corollary}
If $d>1$, then the dimension of $\gch(\mathscr{F}^+)$ is strictly smaller than
the dimension of $\mathscr{C}(\mathscr{F}^+)$.
\end{corollary}

In the case $d=1$, we have five linearly independent one-point functions
by Theorem \ref{theorem-main}.
Recall from Corollary\ref{cor-slf-dim-d1} that $\dim_\C\slf^{\mathscr{T}^+}=5$.
Therefore, by Theorem \ref{theorem3.3.4}, we have:

\begin{theorem}
The space of one-point functions of $\mathscr{T}^+$ is 
$5$-dimensional.
In particular, 
 $\mathscr{C}(\mathscr{T}^+)\cong\gch(\mathscr{T}^+)$.
\end{theorem}

We close this paper by giving the following conjecture.

\begin{conjecture}
$\dim_\C\mathscr{C}(\mathscr{F}^+)=2^{2d-1}+3$
for any positive integer $d$.
\end{conjecture}

\end{document}